\documentclass{article}
    \PassOptionsToPackage{numbers, compress}{natbib}
\usepackage[preprint]{neurips_2025}

\usepackage[utf8]{inputenc} 
\usepackage[T1]{fontenc}    
\usepackage{hyperref}       
\usepackage{url}            
\usepackage{booktabs}       
\usepackage{amsfonts}       
\usepackage{nicefrac}       
\usepackage{microtype}      
\usepackage{xcolor}         
\usepackage{amsthm}         
\usepackage{thmtools}       
\usepackage{amsmath}        
\usepackage{wrapfig}        
\usepackage{mathtools}      
\usepackage{amssymb}

\declaretheorem[numberwithin=section]{theorem}
\declaretheorem[sibling=theorem]{lemma}
\declaretheorem[sibling=theorem]{proposition}

\declaretheorem[sibling=theorem]{remark}
\declaretheorem[sibling=theorem]{definition}
\declaretheorem[sibling=theorem]{assumptions}

\numberwithin{equation}{section}

\title{Discretization Error of Fourier Neural Operators}

\author{%
    Margaret Trautner \\
  Caltech \\
  \texttt{trautner@caltech.edu} \\
  \And
  Samuel Lanthaler \\
  University of Vienna \\
  \texttt{samuel.lanthaler@univie.ac.at} \\
  \And
  Andrew M. Stuart \\
  Caltech \\
  \texttt{astuart@caltech.edu} \\
}

\usepackage{amsfonts}

\usepackage{enumitem}
\usepackage{thm-restate}
\usepackage{todonotes}
\usepackage{subcaption}
\usepackage{datetime}
\usepackage{wrapfig}
\usepackage{type1cm}

\ifpdf
  \DeclareGraphicsExtensions{.eps,.pdf,.png,.jpg}
\else
  \DeclareGraphicsExtensions{.eps}
\fi

\title{Discretization Error of Fourier Neural Operators}

\usepackage{amsopn}

\def\Z{\mathbb{Z}}

\def\Td{\mathbb{T}^d}
\def\R{\mathbb{R}}

\def\cK{\mathcal{K}}
\def\mcE{\mathcal{E}}
\def\mcF{\mathcal{F}}

\def\cA{\mathcal{A}}
\def\cP{\mathcal{P}}
\def\cK{\mathcal{K}}
\def\cQ{\mathcal{Q}}

\def\la{\langle}
\def\ra{\rangle}
\def\R{\mathbb{R}}

\def\XN{X^{(N)}}
\def\mcF{\mathcal{F}}
\def\mcK{\mathcal{K}}
\def\DFT{\mathsf{DFT}}
\def\mcE{\mathcal{E}}
\def\mcEzero{\mathcal{E}^{(0)}}
\def\mcEone{\mathcal{E}^{(1)}}
\def\mcEtwo{\mathcal{E}^{(2)}}
\def\mcEthree{\mathcal{E}^{(3)}}

\def\vhat{\widehat{v}}

\def\ltwoN{\ell^2(n \in [N]^d)}

\def\Zd{\mathbb{Z}^d}
\def\dx{\; \mathsf{d}x}
\def\dy{\; \mathsf{d}y}

\def\psik{e^{2\pi i \langle k, x \rangle}}

\renewcommand{\Omega}{\mathcal{D}}
\newcommand{\defeq}{\coloneqq}

\numberwithin{equation}{section}
\usepackage{color}
\definecolor{darkred}{rgb}{.6,0,0}
\definecolor{darkblue}{rgb}{0,0,.7}
\definecolor{darkgreen}{rgb}{0,.5,0}
\definecolor{darkbrown}{rgb}{0.8,0.4,0.4}

\newcommand{\mt}[1]{{\color{black}{#1}}}

\begin{document}

\maketitle

\begin{abstract}
Operator learning is a variant of machine learning that is designed to 
approximate maps between function spaces from data. 
The Fourier Neural Operator (FNO) is one of the main model architectures used for operator learning. The FNO combines linear and nonlinear operations in physical space with linear operations in Fourier space, leading to a parameterized map acting between function spaces. Although in definition, FNOs are objects in continuous space and perform convolutions on a continuum, their implementation is a discretized object performing computations on a grid, 
allowing efficient implementation via the FFT. Thus, there is a discretization error between the continuum FNO definition and the discretized object used in practice that is separate from other previously analyzed sources of model error. We examine this discretization error here and obtain algebraic rates of convergence in terms of the grid resolution as a function of the input regularity. Numerical experiments that validate the theory and describe model stability are performed. In addition, an algorithm is presented that leverages the discretization error and model error decomposition to optimize computational training time.
\end{abstract}



\section{Introduction}
\label{sec:I}

While most learning architectures are designed to approximate maps between finite-dimensional spaces, operator learning is a method that approximates maps between infinite-dimensional function spaces. These maps appear commonly in scientific machine learning applications such as surrogate
modeling of partial differential equations (PDEs) or model discovery from
data. Fourier Neural Operators (FNOs) are a type of operator learning architecture that parameterize the model directly in function 
space, naturally generalizing deep neural networks (DNNs) \citep{li2021fourier}.
In particular, each 
hidden layer of an FNO assigns a trainable integral kernel that acts on the hidden states by convolution in addition to the usual affine weights and biases of 
a DNN. Taking advantage of the duality between convolution and multiplication under Fourier transforms, these convolutional kernels are represented by Fourier multiplier matrices, whose components are optimized during training
alongside the regular weights and biases acting in physical space.
FNOs have proven to be an effective and popular operator learning method 
in several PDE application areas including weather forecasting 
\citep{pathak2022fourcastnet}, biomedical shape optimization \citep{zhou2024ai}, 
and constitutive modeling \citep{bhattacharya2023learning}.
It is thus of interest to study their theoretical properties.

Although FNOs approximate maps between function spaces, in practice, these functions must be discretized. \mt{In particular, kernel integral operators, including the FNO, perform convolution via an integration that must be computed numerically. The error arising from this difference is called \textit{aliasing error}, and during a forward pass of the FNO, the aliasing error propagates through the subsequent model layers and may be amplified by nonlinearities. Thus, the continuum FNO object differs from the implemented model due to \textit{discretization error.} This may be summarized by the following decomposition: 
 \begin{equation}\label{eqn:err_split}
\Psi^\dagger - \Psi_{FNO}^{N} = 
\underbrace{\left[\Psi^\dagger - \Psi_{FNO} \right]}_{\text{model discrepancy}} 
+ \underbrace{\left[\Psi_{FNO} - \Psi_{FNO}^{N}\right]}_{\text{discretization error}}.
\end{equation}
}Here, $\Psi^\dagger$ is the true map to be approximated by a data-driven model, $\Psi_{\text{FNO}}$ is the continuum FNO map, and $\Psi^N_{\text{FNO}}$ is the discretized version of the FNO. In previous analyses of the universal approximation properties of the FNO, the discretization error component is ignored completely; only the continuum definition of the FNO is used \citep{kovachki2021universal,kovachki2023neural}. While this approach to universal approximation is mathematically sound, it leaves the discretization components of the error unquantified in practice. \mt{Understanding and controlling this discretization error is as important for this model as bounding the model discrepancy error arising from sources such as limited data, optimization, and model capacity.} In this paper, we analyze the discretization error both in theory and experimentally.

Aliasing error depends on the regularity, or smoothness, of the input function in the Sobolev sense; \mt{this is well known in Fourier analysis.} Thus, to bound the error for an entire \mt{FNO implementation}, regularity must be maintained as the state passes through the layers of the network, including the nonlinear activation function. In particular, regularity-preserving properties of compositions of nonlinear functions are required. Bounds of this type are given by \citet{moser1966rapidly} and form a key component of the proofs in this 
work \citep{moser1966rapidly}. Because the smooth GeLU (Gaussian Error Linear Unit) 
\citep{hendrycks2016gaussian} activations preserve regularity, while the 
non-differentiable ReLU activations do not, the analysis in this paper is confined to the former and extends to other smooth activation functions.

\subsection{Contributions}
In this paper, we make the following contributions. 
\begin{enumerate}[label=(C\arabic*),topsep=0.5ex,itemsep=0.5ex,partopsep=1ex,parsep=1ex]
    \item We bound the discretization error that results from \mt{implementing the continuum FNO on a grid.}
    \item We validate this theory concerning the discretization error of the FNO with numerical experiments.
    \item We propose an adaptive subsampling algorithm for faster operator learning training.
\end{enumerate}

In Section \ref{sec:SU} we set up the framework for our theoretical 
results. Section \ref{sec:TR} studies the discretization error of the FNO 
in theory, making contribution (C1). In Section \ref{sec:NE} we present 
numerical experiments that illustrate the theory and propose an algorithm for adaptively refining the discretization during training, making
contributions (C2, C3). We conclude in Section \ref{sec:conc}. \mt{The appendices include a self-contained background on aliasing error as well as additional proofs and technical details.}
\subsection{Related work}
Neural networks have been very successful in approximating solutions of partial differential equations using data. Several approaches are used for such models, including physics-informed neural networks (PINNs), constructive networks, and operator learning models. In the case of PINNs, a standard feed-forward machine learning architecture is trained with a loss function involving a constraint of satisfying the underlying PDE \citep{raissi2019physics}. Another approach to applying machine learning to PDEs is to construct approximating networks from classical PDE-solver methods. For example, in \citep{herrmann2022constructive,Herrmann2020Bayesian,herrmann2024neural, longo2023rham}, ReLU neural networks are shown to replicate polynomial approximations and continuous, piecewise-linear elements used in finite element methods exactly. Both of these two approaches to approximating PDE solution maps require a choice of discretization \mt{within the model} to approximate an infinite-dimensional operator.

Operator learning is a branch of machine learning that aims to approximate maps between function spaces, which include solution maps defined by
partial differential equations (PDEs) \citep{kovachki2023neural}. Several operator learning architectures exist, including DeepONet \citep{lu2021learning}, Fourier Neural Operators (FNO) \citep{li2021fourier}, PCA-Net \citep{bhattacharya2021model}, and random features models \citep{nelsen2021random}. 
Our paper focuses on FNOs, which directly parameterize the model in Fourier space through an integral kernel and allow for changes in discretization in both the input and the output functions, potentially allowing for non-uniform grids \citep{li2023fourier}. 
 In addition, FNO takes advantage of the computational speedup of the FFT to gain additional model capacity with less evaluation time. \mt{A key advantage of the FNO is that it is a \textit{discretization-invariant} operator in the sense that its definition involves no discretization and its implementation can be trivially used on various discretizations with no change of parameter values.}

Error analysis for operator learning begins with establishing universal approximation: results
which guarantee that, for a class of possible maps, a particular choice of model architecture, and a desired maximum error, there exists a parameterization of the model that gives at most that error. Universal approximation results are established for a variety of architectures including ReLU neural networks \citep{cybenko1989approximation}, DeepONet \citep{lanthaler2022error}, FNO 
\citep{kovachki2021universal}, and a general class of neural operators \citep{kovachki2023neural}.  Following universal approximation, model size bounds give a worst-case bound on the model parameter sizes required to achieve a certain error threshold for particular classes of problems. These have been established for FNO \citep{kovachki2021universal,kovachki2024operator}, but the analysis uses only the continuum definition of the FNO and \mt{ignores the fact that in practice the FNO implementation must work with} a discretized version. In this work, we \mt{close the error gap by quantifying and bounding the error that results from discretizing the continuum FNO.}

Perhaps the most conceptually similar work to ours is that of \citet{bartolucci2023neural}, which addresses the fact that discretizations of neural operators deviate from their continuum counterparts. This work introduces  an ``alias-free'' neural operator that bypasses inconsistencies resulting from discretization.
In practice, this research direction has led to operator learning frameworks such as Convolutional Neural Operators (CNO), which are not strictly alias free, but reduce aliasing errors via spatial upsampling \citep{raonic2024convolutional}. These prior works have empirically shown the benefits and importance of carefully controlling discretization errors in operator learning. \mt{Prior work has also examined the effects of changing the number of spectral modes in the implemented FNO and an algorithm to optimize training with variable modes \citep{george2022incremental}. In this work, we also propose an adaptive subsampling algorithm that varies the resolution of the data used in training in a manner designed to minimize training time.}

FNOs remain a widespread neural operator architecture, and an analysis of errors resulting from numerical discretization has so far been missing from the literature. To fill this gap, in this paper we bound the discretization error of FNOs and perform experiments that provide greater insight into the behavior of this error.

\section{Notation}
\label{sec:SU}

Fix integer $m$.
Let $|\cdot|$ denote the Euclidean norm on $\mathbb{R}^m$, $\|\cdot\|$ the $L^2(\Td,\mathbb{R}^m)$ norm, and $\|\cdot\|_{\infty}$ the $L^{\infty}(\Td,\mathbb{R}^m)$ norm. Here, $\Td$ denotes the $d$-dimensional torus, which we identify with $[0,1]^d$ with periodic boundary conditions; we simply write $L^2(\Td)$ or $L^{\infty}(\Td)$ when no confusion will arise. \mt{Let $\|\cdot\|_2$ be the induced matrix 2-norm and $\|\cdot\|_F$ be the Frobenius norm. For a \textit{shallow neural network} $\phi(u) = \sigma(A_1u + b)$ with matrix $A_1$ and vector $b$, we denote by $\|\phi\|^2_{\text{NN}}\defeq \|A_1\|^2_F + |b|^2$.} For nonnegative integer $s$, define the 
Sobolev space $H^s(\Td)=H^s(\Td,\mathbb{R}^m)$ as 
\begin{equation}\label{eqn:Hs}
    H^s(\Td) = \Bigg\{f: \Td \to \mathbb{R}^m \; \Big| \; \sum_{k \in \Zd} ( 1 + |k|^{2s} )| \widehat{f}(k)|^2 < \infty\Bigg\}
\end{equation}
where $\widehat{f}$ denotes the Fourier transform of $f$. Define the semi-norm
\begin{equation*}\label{eqn:seminorm}
    |v|_s^2 : = \int_{\Td} v(-\Delta)^s v \dx
\end{equation*}
for functions $v:\Td \to \mathbb{R}^m$. It is useful to consider the following equivalent definition of the space $H^s(\Td)$ for integer $s > d/2$ in terms of this seminorm:
\begin{align*}
    H^s(\Td) &= \{f: \Td \to \R^m \; | \; \|f\|_{H^s}< \infty\} \\
    \|f\|_{H^s} & = \left((2\pi)^{-2s} |f|_s^2 + \|f\|^2\right)^{1/2}.
\end{align*}
We say an element $f \in H^{s-}$ if $f \in H^{s-\epsilon}$ for any $\epsilon >0$. Further, let  $X^{(N)}$ denote the $d$-dimensional grid $\frac{1}{N}[N]^d$ where \[[N]^d := \{n \in \mathbb{Z}_{\geq 0}^d\; | \; n_i < N, \; i \in \{1, \dots, d\}\}.\] Here, $n_i$ is the $i$th entry of vector $n$. We assume $N>1$ throughout this work. We also introduce the following symmetric index set for the Fourier coefficients: $[[N]]^d = [[N]]\times \dots \times [[N]]$, where 
\[
[[N]] := 
\begin{cases}
\{-K,\dots, K\}, & (N=2K+1 \text{ is odd}), \\
\{-K,\dots, K-1\}, & (N=2K \text{ is even}).
\end{cases}
\]
Irrespective of whether $N$ is odd or even, $[[N]]^d$ contains $N^d$ elements.
For functions $u: \Td \to \R^m$, we abuse notation slightly and use $\|u\|_{\ell^2(n\in[N]^d)}$ to indicate the quantity,
\begin{equation*}
\|u\|_{\ell^2(n\in[N]^d)} := \Big(\sum_{n \in [N]^d}|u(x_n)|^2\Big)^{1/2}.
\end{equation*}
This is a norm for the vector found by evaluating $u$ at grid points. Note that for $x_n = \frac{1}{N}n$ where $n \in [N]^d$, it holds that $x_n \in \Td$, and if $u \in L^2(\Td)$ is Riemann integrable,
\begin{equation}
\label{eq:2L2}
    \lim_{N \to \infty} \frac{1}{N^{d/2}}\|u\|_{\ell^2(n \in [N]^d)} = \|u\|_{L^2(\Td)}.
\end{equation}

Finally, we define the FNO. We remark that this constitutes the standard definition of the FNO with the exception that we
ask for smooth activation functions. At a high level, the FNO is a composition of layers, where the first and final layers are lifting and projection maps, and the internal layers are an activation function acting on the sum of an affine term, a nonlocal integral term, and a bias term. \mt{We emphasize that this FNO definition \textit{does not involve a discretization}; it is a map between function spaces on a continuum.}
\begin{definition}[\bf{Fourier Neural Operator}] Let $\mathcal{A}$ and $\mathcal{U}$ be two Banach spaces of real vector-valued functions
over domain $\Td$. Assume input functions $a \in \mathcal{A}$ are $\R^{d_a}$-valued while the output functions $u \in \mathcal{U}$ are $\R^{d_u}$-valued. The neural operator architecture $\Psi_{\theta}:\mathcal{A} \to \mathcal{U}$ is
\begin{align*}
    \Psi_{\theta} &= \mathcal{Q} \circ \mathsf{L}_{T-1}\circ \dots \circ \mathsf{L}_0\circ \mathcal{P},\\
    v_{t+1} &= \mathsf{L}_t v_t  =\sigma_t(W_{t}v_t + \mathcal{K}_{t}v_t + b_{t}), \quad t=0,1, \dots, T-1,
\end{align*}
with \(v_0 = \mathcal{P}(a)\). Here, $\mathcal{P}: \R^{d_a} \to \R^{d_{0}}$ and $\mathcal{Q}:\R^{d_{T}}\to\R^{d_u}$ are \mt{shallow neural networks with globally Lipschitz, $C^{\infty}$ activations $\sigma_p$ and $\sigma_q$,} and the $\sigma_t$ are fixed nonlinear activation functions acting locally as maps $\R^{d_{t+1}} \to \R^{d_{t+1}}$ in each layer.  $\mathcal{P}$, $\mathcal{Q}$ and the $\sigma_t$ are viewed
as operators acting pointwise, or pointwise almost everywhere, over the domain $\Td$, $W_t \in \R^{d_{t+1}\times d_{t}}$ are matrices, $\mathcal{K}_t: \{v_t: \Td \to \R^{d_{t}}\} \to \{{v_{t+1}}:\Td \to \R^{d_{t+1}}\}$ are integral kernel operators and $b_t$ is a bias term. The activation functions $\sigma_t$ are restricted to the set of globally Lipschitz, non-polynomial, $C^{\infty}$ functions. The integral kernel operators $\mathcal{K}_t$ are parameterized in the Fourier domain in the following manner.  Let \(i = \sqrt{-1}\) denote the imaginary unit.
Then, for each $t$, the kernel operator $\mathcal{K}_t$ is parameterized by 
\begin{equation}\label{eqn:kernel}
    (\mathcal{K}_tv_t)(x) = \sum_{k \in [[K]]^d}\sum_{j=1}^{d_{t}} (P^{(k)}_t)_{j}\langle e^{2\pi i \langle k,\cdot \rangle},(v_t)_j\rangle_{L^2(\mathbb{T}^d;\mathbb{C})} \psik \in\R^{d_{t+1}}\,.
\end{equation}
Here, each $P_t^{(k)} \in \mathbb{C}^{{d_{t+1}} \times d_{t}}$ constitutes the learnable parameters of the integral operator, with $(P_t)^{(k)}_j$ the $j$th column, and $K \in \mathbb{Z}^+$ is a mode truncation parameter. \mt{$\mathcal{K}_t$ is well-defined for $v_t \in L^2(\Td)$.} We denote by $\theta$ the collection of parameters that specify $\Psi_{\theta}$, which include the weights $W_t$, biases $b_t$, kernel weights $P_t$, and the parameters describing $\mathcal{P}$ and $\mathcal{Q}$.
\label{def:FNO}
\end{definition}

In the error analysis in the following section, we are interested in the discrepancy between taking the inner product in equation \eqref{eqn:kernel} on a grid instead of on a continuum -- the errors due to \emph{aliasing}. \mt{The above continuum definition is assumed in learning theory analysis of the FNO, but in practice, a discretized approximation is used.} We consider the other parameters, including the mode count $K$, to be fixed
and intrinsic to the FNO model considered, irrespective of which grid it is approximated on.

\section{Main results} 
\label{sec:TR}

 \mt{Let $\mathcal{A}$ and $\mathcal{U}$ be the input and output Banach spaces in the FNO definition \ref{def:FNO}, and let the FNO model hyperparameters be fixed. Given a setting of the trainable FNO parameters $\theta$, let $\Psi_{\text{FNO}}:\mathcal{A} \mapsto \mathcal{U}$ be the FNO obtained using the definition. \textit{This definition does not involve a discretization.} Thus, any implementation of the FNO with the same hyperparameters and trainable $\theta$ must be another map, denoted $\Psi^{N}_{\text{FNO}}: \mathcal{A} \mapsto \mathcal{U}$, that evaluates the $L^2$ inner product in equation \eqref{eqn:kernel} numerically on grid points $X^{(N)}$ rather than at every point $x \in \Td$ as $\Psi_{\text{FNO}}$ does. In particular, $\Psi^N_{FNO}$ exchanges the operator $\cK_t$ defined in \eqref{eqn:kernel} for $\cK_t^N$ such that \begin{equation}\label{eqn:KN}
     (\cK^N_t v^N_t)(x_n) = \sum_{k \in [[K]]^d}\sum_{j=1}^{d_{t}} (P^{(k)}_t)_{j}\DFT((v_t^N)_j)(k) e^{2\pi i \langle k, x_n \rangle} \in\R^{d_{t+1}},
 \end{equation}
 where $\DFT$ is the discrete Fourier transform. We refer to the output of each internal layer $\mathsf{L}$ as a \textit{state value}. Starting from the same input $a \in \mathcal{A}$, $\Psi_{\text{FNO}}(a)$ and $\Psi_{\text{FNO}}^N(a)$ will have different state values, denoted $v_t$ and $v_t^N$, respectively, as outputs of internal layer $\mathsf{L}_{t-1}$ for $t>0$ despite having the exact same model parameters. This difference is important because in proofs concerning the FNO, only $\Psi_{\text{FNO}}$ is considered, but $\Psi_{\text{FNO}}$ is an \textit{unimplementable} object in practice. If $\Psi^\dagger$ is the map of interest to be approximated using an FNO model, the overall approximation error of the implemented $\Psi_{\text{FNO}}^N$ can be split into a contribution due to the numerical discretization and another contribution due to model discrepancy, as shown in \eqref{eqn:err_split}.
 Theorem \ref{thm-fno:main} bounds the discretization error component. The result takes into account both the initial errors that occur with each approximated inner product and their magnified effects as they propagate through the layers of the model. Despite this nonlinear propagation, we show that the approximate $L^2$ norm of the error after any number of layers decreases like $N^{-s}$, where $s$ describes the Sobolev regularity of the input.}


To prove Theorem \ref{thm-fno:main}, we assume the following.

\begin{assumptions} \label{asst:main} For a fixed FNO with $T$ layers:
\vspace{-0.3cm}
    \begin{enumerate}[label=(A\arabic*),topsep=1.67ex,itemsep=0.1ex,partopsep=1ex,parsep=0.5ex]
        \item There exists some $B \geq 1$ such that $\sigma_t,\sigma_p,$ and $\sigma_q$ possess continuous derivatives up to order $s$ which are bounded by $B$. \label{ass:A1}
        \item Input set $\cA \subset H^s(\Td)$. 
        \item $1 \leq K < \frac{N}{2}$. 
        \item $s > \frac{d}{2}$. 
        \item There exists some $M \geq 1$ such that FNO parameters $P_t$, $W_t$, and $b_t$ are each bounded above by $M$ in the following norms: $\|P_t\|_F\leq M$, $\|W_t\|_2 \leq M$, and $|b_t|\leq M$ for all $t \in [0, \dots, T-1]$. \mt{Furthermore, $\cP$ and $\cQ$ are bounded and smooth with $\|\cP\|_{\text{NN}} \leq M$, and $\|\cQ\|_{\text{NN}} \leq M$. }
    \end{enumerate}
\end{assumptions}

 The above assumptions are easily satisfied in practice. We investigate the consequences of violating Assumption \ref{ass:A1} in the numerical experiments. The main result is the following theorem concerning the behavior of the error with respect to the size of the discretization. 

\begin{restatable}{theorem}{thmmain}\label{thm-fno:main} Let Assumptions \ref{asst:main} hold. \mt{Let $\cA_c$ be a compact set in $\cA$. Let $v_t(a)\defeq \mathsf{L}_{t-1}\circ\dots \circ \mathsf{L}_0\circ \cP(a)$ with $\cP$ and each $\mathsf{L}$ as defined in Definition \ref{def:FNO}.}  Similarly, let $v_t^N(a)\defeq \mathsf{L}^N_{t-1} \circ\dots \circ \mathsf{L}^N_0\circ \cP(a)$ where $\mathsf{L}^N_jv_j^N = \sigma_j(W_jv^N_j + \cK^N_jv^N_j+b_j)$ for $\cK_j^N$ defined in \eqref{eqn:KN} for each $0 \leq j \leq t-1$. Then
\begin{align}\label{eqn:decay}
        \sup_{a\in \cA_c} \frac{1}{N^{d/2}}\|v_t(a) - v^N_t(a)\|_{\ltwoN}
        & \leq 
        CN^{-s}
\end{align} 
where the constant C depends on $B,M,d,s,t,$ and $\cA_c$.
\end{restatable}
The proof and exact form of the constant $C$ in the above theorem are detailed in the supplementary materials. 

We can also state the following variant of Theorem \ref{thm-fno:main}, which shows that the same convergence rate is obtained at the continuous level, when $v^{N}_t(x_n)$ is replaced by a trigonometric polynomial interpolant:
\begin{restatable}{theorem}{thmmainp}\label{thm:main'}
Let $p^N_t(a)(x) = \sum_{k\in [[N]]^d} \DFT(v^N_t(a))(k) e^{2\pi i \langle k, x\rangle}$ denote the interpolating trigonometric polynomial of $\{v^N_t(a)(x_n)\}_{n\in [N]^d}$. Let the assumptions of Theorem \ref{thm-fno:main} hold. Then,
\begin{align}\label{eqn:decay'}
       \sup_{a\in\cA_c}\Vert v_t(a)- p^N_t(a)\Vert_{L^2(\Td)} &\le C' N^{-s}.
\end{align}
Here, $C'$ depends on $B,M,d,s,t,$ and $\cA$. 
\end{restatable}
\mt{
\begin{remark}
    A consequence of Theorem \ref{thm:main'} is that $\sup_{a\in\cA_c}\|\Psi_{\text{FNO}}(a) - \Psi^N_{\text{FNO}}(a)\|_{L^2(\Td)}$ also has a convergence rate of $N^{-s}$ since under Assumptions \ref{asst:main} $\cQ$ is Lipschitz and $\cQ,\cP$ preserve regularity.
\end{remark}
}



The exact form of the constant $C'$ may be found in the proof in the appendix. A key element in the proof of Theorem \ref{thm-fno:main} is to provide a bound on the Sobolev norm of the ground truth state $\|v_t\|_{H^s}$ at each layer. The following lemma accomplishes this for a single layer. The proof may be found in the appendix.
\begin{restatable}{lemmma}{lemHs} \label{lem:HsBound}
    Under Assumptions \ref{asst:main}, the following bounds hold: 
    \begin{itemize}[topsep=1.67ex,itemsep=0.5ex,partopsep=1ex,parsep=1ex]
        \item $\|v_{t+1}\|_{\infty} \leq \sigma^* + BM(1+\|v_t\|_{\infty} + K^{d/2}\|v_t\|_{L^2(\Td)}).$
        \item $|v_{t+1}|_s \leq BcM^sK^{ds/2}(1+\|v_t\|_{\infty})^s(1+|v_t|_s)$
    \end{itemize}
    for some constant $c$ dependent on $d$ and $s$, where $\sigma^* \defeq \max\{\max_{0 \leq t \leq T-1} \sigma_t(0),1\} $. 
\end{restatable}

The result of Theorem \ref{thm-fno:main} guarantees that the discretization error 
converges as grid resolution increases. The algebraic decay rate in
a discrete $L^2$ norm is determined by the regularity of the input; this
in turn builds on Lemma \ref{lem:HsBound} which ensures that the regularity 
of the state is preserved through each layer of the FNO.

\section{Numerical experiments}
\label{sec:NE}

In this section we present and discuss results from numerical experiments that empirically validate the results of Theorem \ref{thm-fno:main}. The $L^2$ error at each layer decreases like $N^{-s}$ where $s$ governs the input regularity and $N$ is the discretization used to perform convolutions in the FNO implementation. For each FNO model in this section, we use a computation 
of a discrete FNO on a high resolution grid as the ``ground truth;'' \mt{this is standard practice in numerical analysis when the true solution is unobtainable.}
We compare states at each layer resulting from inputs of lower resolution with the state resulting from the ground truth. To obtain evaluations of $v_{\ell}$ at higher discretizations than $N$, the inverse Fourier transform operation is interpolated to additional grid points using trigonometric polynomial interpolation; Theorem \ref{thm:main'} justifies this practice.

We perform experiments for inputs of varying regularity by generating Gaussian random field (GRF) inputs with prescribed smoothness $H^{s-}$ for $s \in \{0.5, 1, 1.5, 2\}$. The GRF inputs are discretized for values of $N \in \{32, 64, 128, 256, 512, 1024, 2048\}$ and $d=2$. Grid size $2048$ is used as the ground truth, and the relative error at layer $\ell$ for $v_{\ell}^N$ compared with the truth $v_{\ell}$ is computed with 
\begin{equation*}
    \label{eqn:rel_err}
    \text{Relative Error = } \frac{\|v_{\ell} - v^N_{\ell}\|_{\ell^2(n \in [2048]^d)}}{\|v_{\ell}\|_{\ell^2(n \in [2048]^d)}}.
\end{equation*}

Finally, in FNO training, it is common practice to append positional information about the domain at each evaluation point in the form of Euclidean grid points; i.e. $(x_1,x_2) \in [0,1]^2$ for two dimensions. However, this grid information is not periodic, and an alternative is to append periodic grid information; i.e.$(\sin(x_1), \cos(x_1), \sin(x_2), \cos(x_2))$ for two dimensions. In these experiments, we also compare the error of models with these two different positional encodings.

In Subsection \ref{ssec:random} we discuss experiments on FNOs with random weights, and in Subsection \ref{ssec:trained} we discuss experiments on trained FNOs.\mt{ In the random weights experiments, we present a few interesting experimental findings; namely, using ReLU activations or non-periodic position encodings negatively affects the discretization error decay as the theory predicts. In the trained network experiments, we explore the example of learning a gradient map to show that the model cannot learn an output with less regularity than the input.} Finally, in Subsection \ref{ssec:adaptive}, we propose an application of discretization subsampling to speed up operator learning training by leveraging adaptive grid sizes.

\subsection{Experiments with random weights}\label{ssec:random}

In this subsection, we consider FNOs with random weights and study their discretization error and model stability with respect to perturbations of the inputs. All models are defined in spatial dimension $d = 2$, with $K = 12$ modes in each dimension, a width of $64$, and $5$ layers. 

The default model has randomly initialized iid $\mathcal{U}(-\frac{1}{\sqrt{d_t}}, \frac{1}{\sqrt{d_t}})$ weights (uniformly distributed) for the affine and bias terms, where $d_t$ is the layer width, and iid $\mathcal{U}(0,\frac{1}{d_t^2})$ spectral weights. Initializing the weights this way is the standard default for FNO. This model uses the GeLU activation function standard in FNO. Next we examine the use of ReLU activation instead of GeLU. Finally, we investigate non-periodic positional encoding. 
\begin{figure}[ht!]
    \begin{center}
    \vspace{-10pt}
    \includegraphics[width = \linewidth]{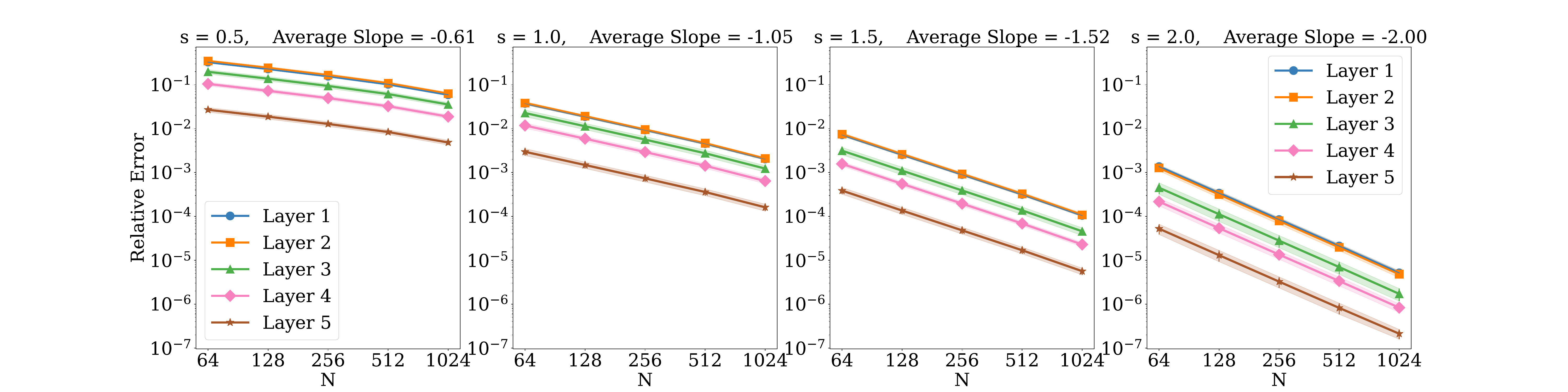}
    \caption{Relative error versus $N$ and $s$ for an FNO with default weight initialization.}
    \label{fig:err_vs_N}
    \end{center}
    \vspace{-20pt}
\end{figure}

\begin{figure}[hb!]
    \centering
    \vspace{-10pt}
    \includegraphics[width = \linewidth]{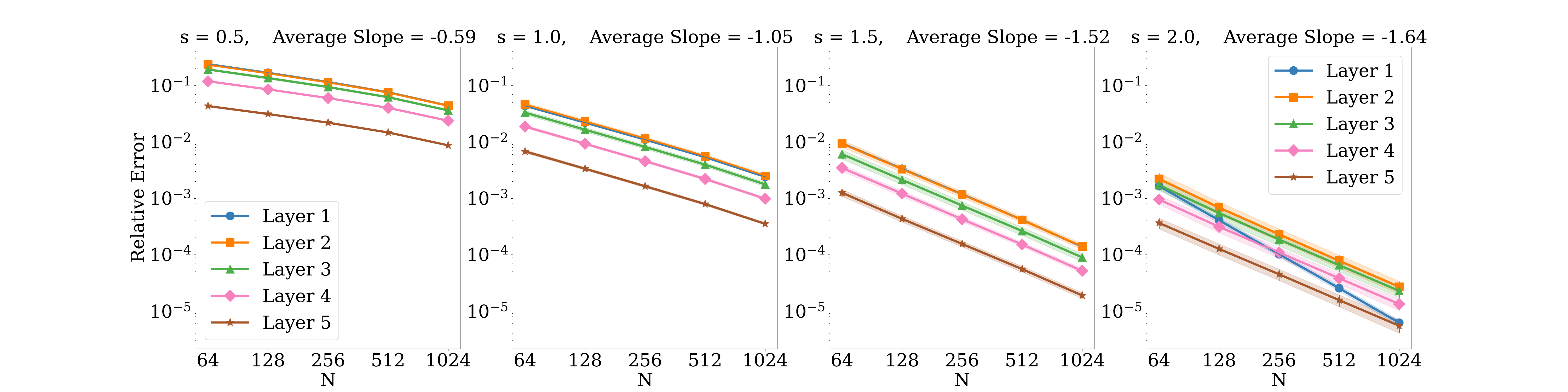}
    \caption{Relative error versus $N$ and $s$ for a default FNO with a ReLU activation.}
    \label{fig:err_vs_N_relu}
    \vspace{-5pt}
\end{figure}

\begin{figure}[ht!]
    \centering
    \vspace{-10pt}
    \includegraphics[width = \linewidth]{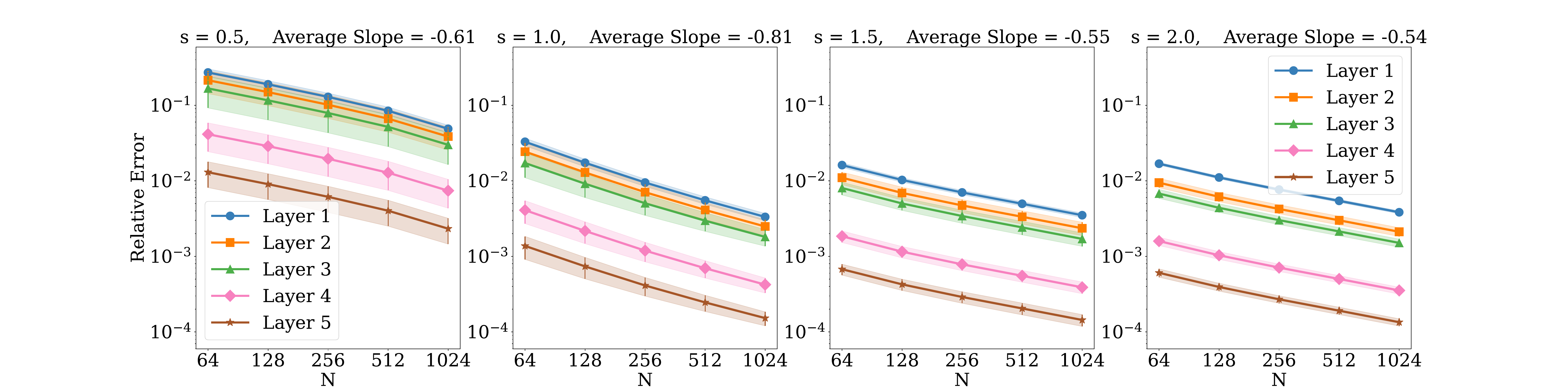}
    \caption{Relative error versus $N$ and $s$ for a default FNO with non-periodic position encoding appended to the input.}
    \label{fig:err_vs_N_xygrid}
    \vspace{-10pt}
\end{figure}

\paragraph{Discretization error for random weights models}
The relative error of the state at each layer versus the discretization for inputs of varying regularity may be seen for the default model, the ReLU model, and the non-periodic position encoding model in Figures \ref{fig:err_vs_N},\ref{fig:err_vs_N_relu}, and \ref{fig:err_vs_N_xygrid} respectively. In these figures, from left to right, $s \in \{0.5, 1, 1.5, 2\}$ where $v_0 \in H^{s-}$. The uncertainty shading indicates two standard deviations from the mean over five inputs to the FNO.

As can be seen in Figure \ref{fig:err_vs_N} for the model with the default weight initialization, the empirical behavior of the error matches the behavior expected from Theorem \ref{thm-fno:main}. One question that arises from Figure \ref{fig:err_vs_N} is why the error decreases as the number of layers increases; this is an effect of the magnitude of the weights. When the model weights are multiplied by $10$, then the error begins to increase with the number of layers. A figure illustrating this phenomenon may be found in the supplementary materials, where additional weight initializations are explored as well.

The results shown in Figure \ref{fig:err_vs_N_relu} justify the use of the GeLU activation function, which belongs to $C^{\infty}$, over the ReLU 
activation function, which is only Lipschitz. The figure shows that the
benefit of having sufficiently smooth inputs is negated by the ReLU activation: the error decay is limited. Note that this effect does not occur for the first layer since at that point ReLU has been applied once, and the Fourier transform is not applied to the output of an activation function until the second layer. \mt{Additionally, we do not observe the effect of ReLU until the input has regularity greater than $s=1.5$ since the ReLU activation function has regularity of $s = 1.5$}

A similar effect to the ReLU model occurs when non-periodic positional encoding information is appended to the input, \mt{as is standard in practical FNO usage;} see Figure \ref{fig:err_vs_N_xygrid}. Since this grid data has a jump discontinuity across the boundary of $[0,1]^d$, it has regularity of $s = 0.5$, so the convergence rate never achieves $N^{-1}$. \mt{These results suggest caution when using positional encoding information with smooth input data; periodic positional encodings may be preferred. }

\subsection{Experiments with trained networks}\label{ssec:trained}

In this subsection, we consider two different maps and train FNOs on data from each map. The first map is a PDE solution map in two dimensions whose solution is at least as regular as the input function. The second map is 
a simple gradient, but in this setting the output function of the gradient is naturally less regular by one Sobolev smoothness exponent than that of the input function. In both experiments, periodic positional encoding information is appended to the inputs. 

\paragraph{Example 1: PDE solution model}
\begin{figure}
\centering
\begin{subfigure}{.5\textwidth}
  \centering
  \includegraphics[width=\linewidth]{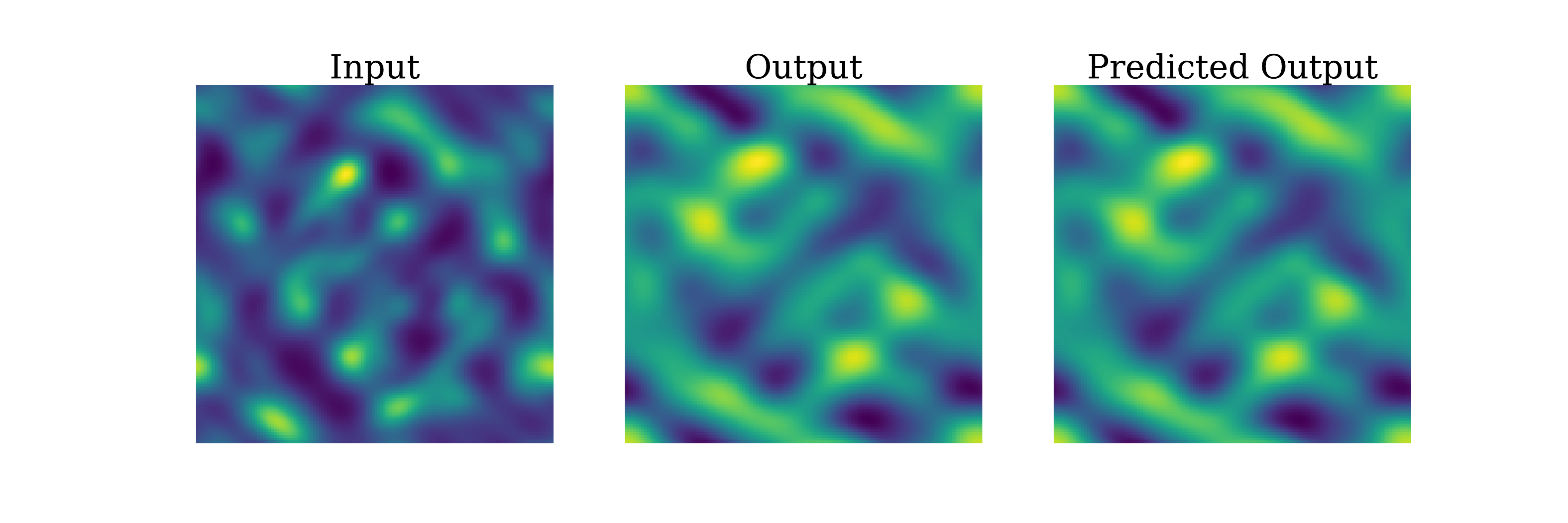}
  \caption{Data for the PDE Solution FNO.}
  \label{fig:data_pde}
\end{subfigure}%
\begin{subfigure}{.5\textwidth}
  \centering
  \includegraphics[width=\linewidth]{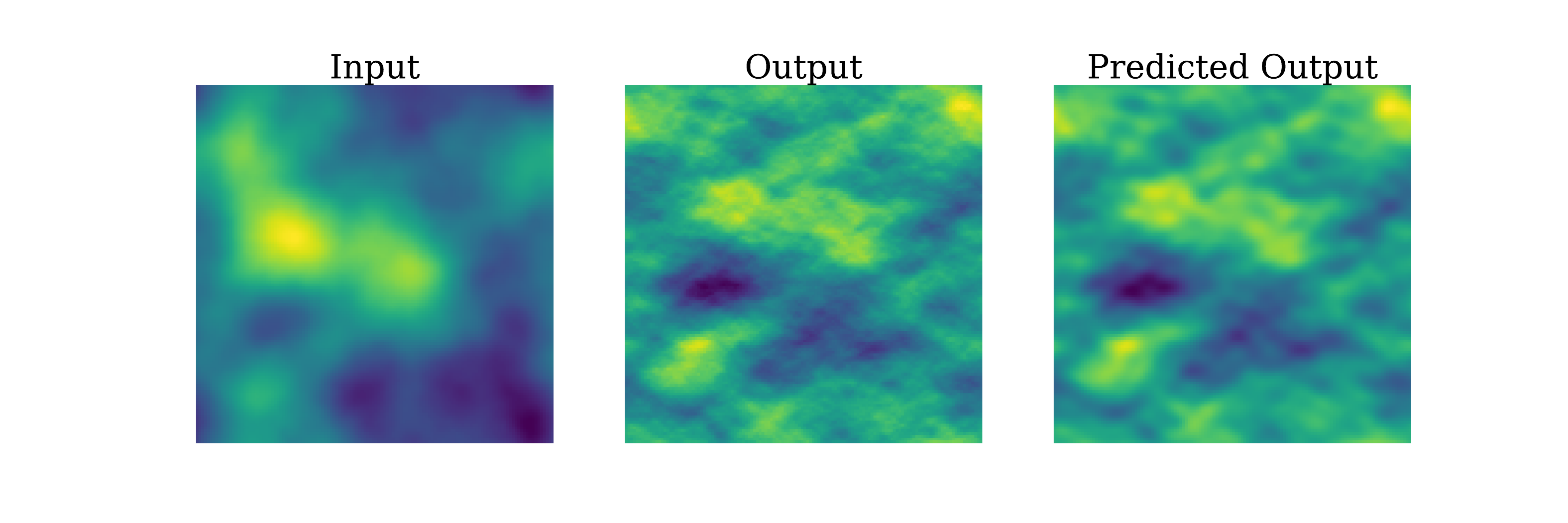}
  \caption{Data for the Gradient FNO.}
  \label{fig:data_grad}
\end{subfigure}
\caption{Visualization of the input and output data for the trained model examples.}
\label{fig:data}
\end{figure}

\begin{figure}
    \centering
    \includegraphics[width = \textwidth]{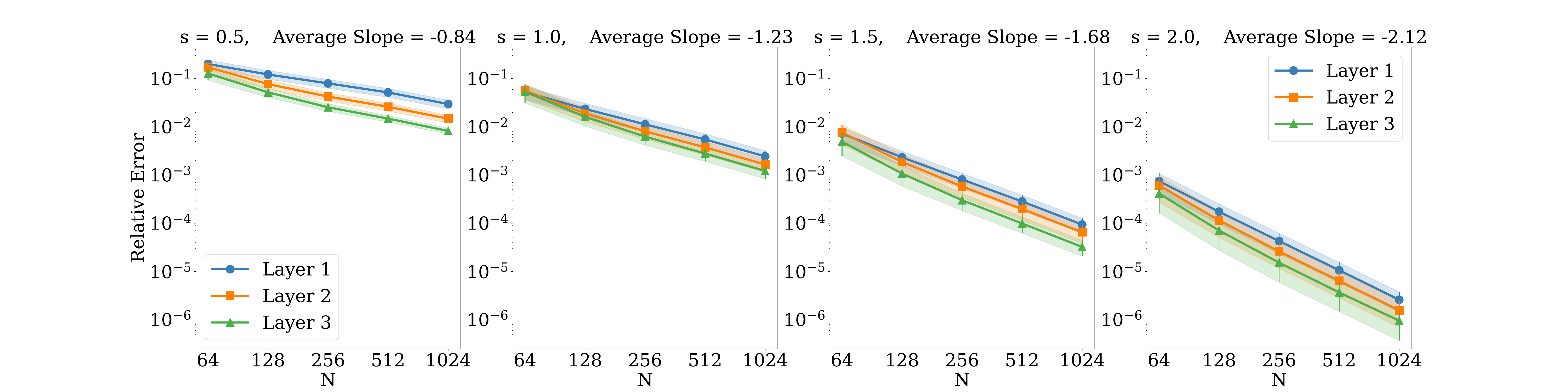}
    \caption{Error versus discretization for inputs of varying regularity for the FNO trained on data corresponding to a PDE solution.}
    \label{fig:smooth_per}
    \vspace{-10pt}
\end{figure}

\begin{figure}
    \centering
    \vspace{-20pt}
    \includegraphics[width = \textwidth]{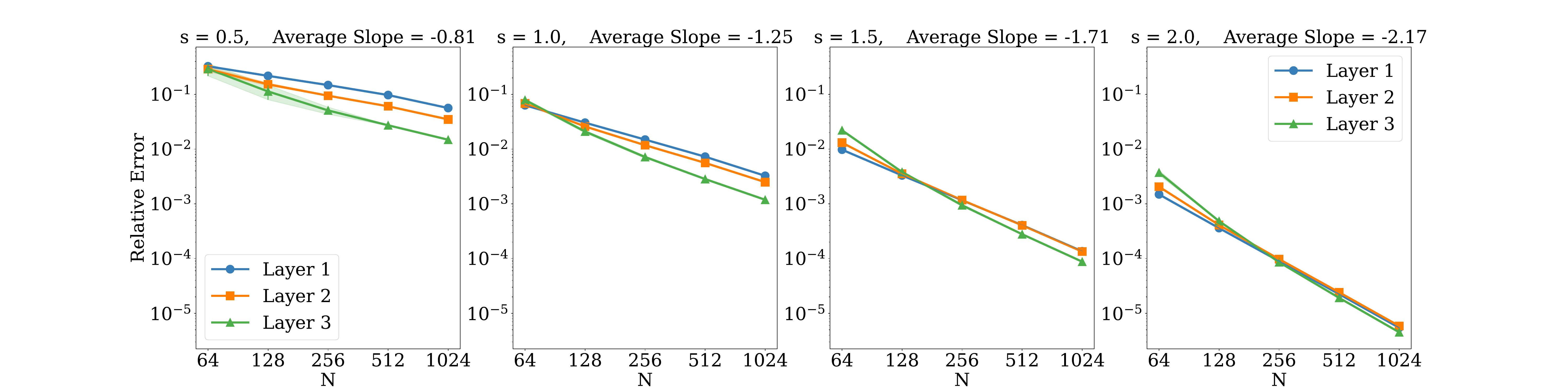}
    \caption{Error versus discretization for inputs of varying regularity for the FNO trained on data corresponding to a gradient map.}
    \vspace{-15pt}
    \label{fig:err_grad}
\end{figure}

In this example, we train an FNO to approximate the solution map of the following PDE:
\begin{align}\label{eqn:cell}
    \nabla\cdot(\nabla \chi A) & = \nabla \cdot A, \quad y \in \mathbb{T}^2 \\ \chi \text{ is } &1-\text{periodic, }\int_{\mathbb{T}^2}\chi \dy = 0.
\end{align}
Here, the input $A: \mathbb{T}^2 \mapsto \R^{2 \times 2}$ is symmetric positive definite at every point in the domain $\mathbb{T}^2$ and is bounded and coercive. For the output data we take the first component of $\chi: \mathbb{T}^2 \mapsto \R^2$.  In our experiments the model is trained to $< 5\%$ relative 
$L^2$ test error. A visualization of the data is in Figure \ref{fig:data_pde}.

The error versus discretization analysis can be seen in Figure \ref{fig:smooth_per}. The error decreases slightly faster than predicted by the
theory; a potential explanation is that the trained model itself 
has a smoothing effect that is not exploited in our analysis.

\paragraph{Example 2: gradient map}

In the final example, we train an FNO to approximate a simple gradient map $u \mapsto \nabla u.$

The training data consists of iid Gaussian random field inputs with regularity $s=2$. Since a gradient reduces regularity, we expect the model outputs to approximate functions with regularity $s = 1$, which is at odds with the smoothness-preserving properties of the FNO described by theory. 

The error versus discretization for inputs of various smoothness is shown in Figure \ref{fig:err_grad}. The error decreases according the the smoothness of the input despite the smoothness-decreasing properties of the data. Indeed, the model does produce more regular predicted outputs than the true gradient, as can be seen in Figure \ref{fig:data_grad} where the predicted output is visibly smoother than the true output. 

\subsection{Speeding up training via adaptive subsampling} \label{ssec:adaptive}

\begin{wrapfigure}{r}{0.35\linewidth}
  \centering
  \vspace{-20pt}
  \includegraphics[width=\linewidth]{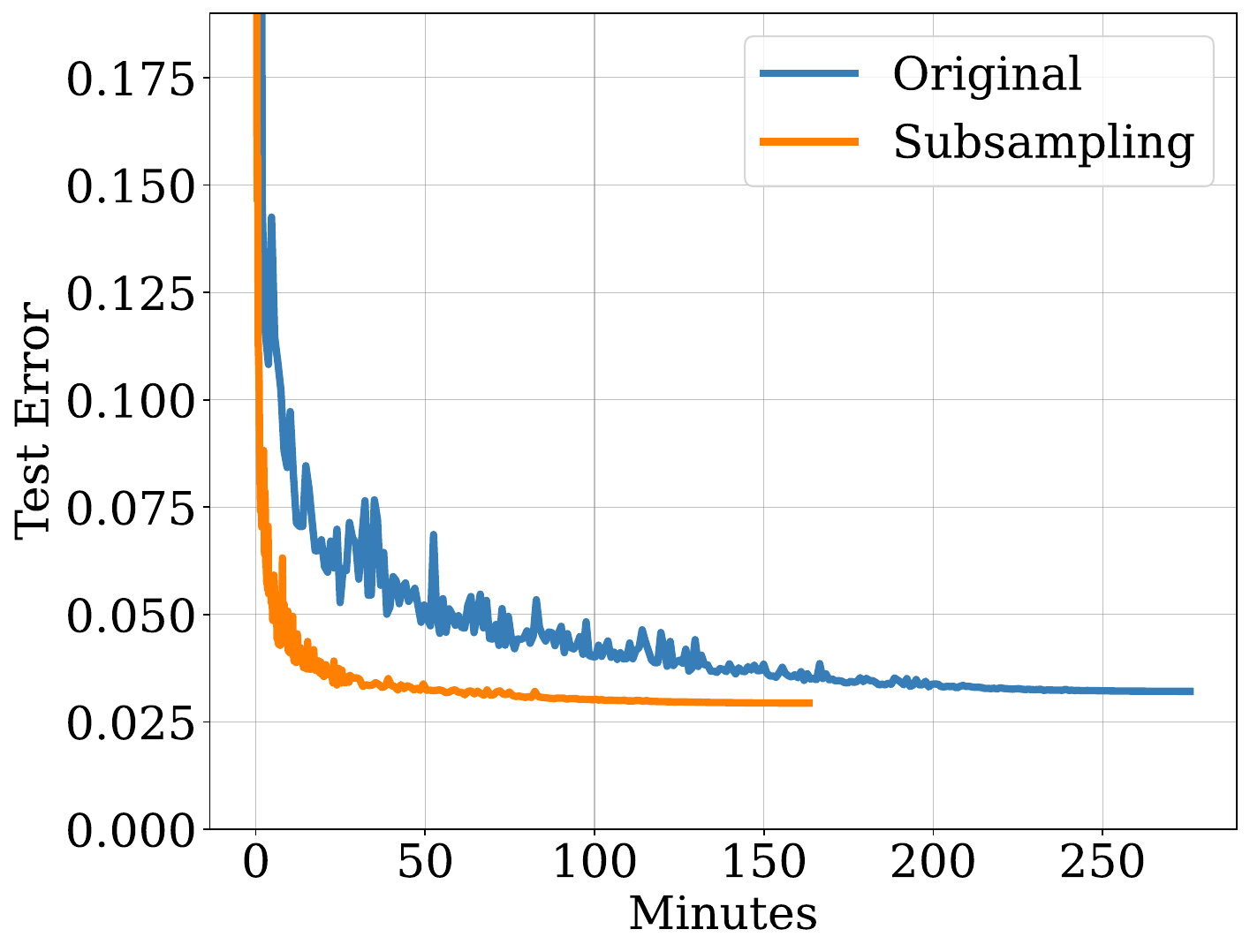}
  \caption{Adaptive grid refinement leads to greater training efficiency.}
  \vspace{-10pt}
  \label{fig:subsampling}
\end{wrapfigure}

The fact that the FNO architecture and its parametrization are independent of the numerical discretization allows for increased flexibility. Specifically, it is possible to adaptively choose an optimal discretization for a given objective. Furthermore, the error decomposition \eqref{eqn:err_split} invites an exploitation to optimize computational training time.

The idea of the proposed approach is that, during training, it is not necessary to compute outputs to a numerical accuracy that is substantially better than the model discrepancy. This suggests an adaptive choice of the numerical discretization, where we use a coarse grid during the early phase of training and refine the grid later. In practice, we realize this idea by introducing a subsampling scheduler. The subsampling scheduler tracks a validation error on held out data and adaptively changes the numerical resolution via suitable subsampling of the training data. Starting from a coarse resolution, we iteratively double the grid size once the validation error plateaus. 

We train FNO for the elliptic PDE \eqref{eqn:cell} with and without the subsampling scheduler; details are contained in the supplementary materials. Compared to training without subsampling, training with a subsampling scheduler requires the same number of forward and backward passes over the network for the training and test set, plus an additional overhead due to the validation set. Since we are mainly interested in the training time, our choice of adding validation samples, rather than performing a training/validation split of the original training samples, ensures that computational timings are not skewed in favor of subsampling. Over the course of training, we iterate through the following grid sizes: 32x32, 64x64, and 128x128. Our criterion for a plateau is that the validation error has not improved for 40 training epochs. The results of training with and without subsampling scheduler for the PDE solution model \eqref{eqn:cell} are shown in Figure \ref{fig:subsampling}. We observe that training time can be substantially reduced with subsampling. \mt{The presentation of the algorithm here is preliminary and serves to emphasize the importance of considering discretization error in the model. The success of this simple example } points to the potential benefits of developing adaptive numerical methods for model evaluation within operator learning.

\section{Conclusions}
\label{sec:conc}

In this paper, we analyze the \mt{discretization error that results from implementation of
Fourier Neural Operators (FNOs).} We bound the $L^2$ norm of the error in 
Theorem \ref{thm-fno:main}, proving an upper bound that decreases 
asymptotically as $N^{-s}$, where $N$ is the discretization in each dimension, 
and $s$ is the input regularity. We show empirically that FNOs with random
weights chosen as the default FNO weights for training behave almost exactly as the theory predicts. Furthermore, our theory and experiments
justify the use of the GeLU activation function in FNO over ReLU, as the former preserves regularity. Additional analyses on trained models show that the error behaves less predictably in relation to our theory
in the low-discretization regime. \mt{Finally, we use the decomposition of model error and discretization error} to propose an adaptive subsampling algorithm for decreasing training time with operator learning. As FNOs become a more common tool in scientific machine learning, understanding the various sources of error is critical. By bounding 
FNO discretization error and demonstrating its behavior in numerical experiments, we understand its effect on learning and the potential to minimize computational costs by an adaptive choice of numerical resolution.

\begin{ack}SL is supported by Postdoc.Mobility grant P500PT-206737 from the Swiss National Science Foundation. The work of AMS is supported by a Department of Defense Vannevar Bush Faculty Fellowship, and by the SciAI Center, funded by the Office of Naval Research (ONR), under Grant Number N00014-23-1-2729. MT is supported by the Department of Energy Computational Science Graduate Fellowship under award number DE-SC00211.  All code and data for this work is available at \url{https://github.com/mtrautner/BoundFNO}. The authors are grateful to Nicholas Nelsen for helpful discussions on FNO implementation. The computations presented here were conducted in the Resnick High Performance Computing Center, a facility supported by Resnick Sustainability Institute at the California Institute of Technology.
\end{ack}

\bibliography{references}
\bibliographystyle{abbrvnat}


\newpage
\appendix
\addcontentsline{toc}{section}{Appendices}
\section*{Appendices}

\section{Trigonometric interpolation and aliasing}\label{apdx:aliasing}

In this section, we present a self-contained analysis of aliasing errors 
for $v \in H^s(\Td)$. \mt{These results are straightforward and well known in numerical analysis, but we give a clear exposition here as background as it is difficult to find a succinct and widely-available reference.} The primary goal is to state and prove 
Proposition \ref{prop:vp}, which controls the difference between a function defined over $\Td$ and the trigonometric interpolation of a function defined on a grid. In the following, $N \in \Z_{>0}$.
We recall that $X^{(N)}$ is a set of equidistant grid points on the torus $\Td$,
\[
X^{(N)} = \{ x_n \in \Td \, |\, x = n/N, \; n \in [N]^d \}.
\]
We note that the discrete Fourier transform gives rise to a natural correspondence between grid values and Fourier modes,
\begin{equation}
\label{eq:DFTd}
\{v(x_n)\}_{n \in [N]^d} \; \leftrightarrow \; \{ \widehat{v}_k \}_{k\in [[N]]^d},
\end{equation}
where 
\begin{equation}\label{eqn:DFT}
\widehat{v}_k = \frac{1}{N^d} \sum_{n \in [N]^d} v(x_n) e^{-2\pi i \langle k, x_n\rangle} =: \DFT(v)(k).
\end{equation}
We begin with the following observation:

\begin{lemma}
Let $N$ be given. Then, 
\begin{align}
\frac1{N^d} \sum_{k \in [[N]]^d} e^{2\pi i\langle k, x_m - x_n \rangle} = \delta_{mn},  \quad \forall \, m,n \in [N]^d, 
\label{eq:t1}
\\
\frac1{N^d} \sum_{n \in [N]^d} e^{2\pi i\langle k-k', x_n \rangle} = \delta_{kk'}, \quad \forall \, k,k' \in [[N]]^d.
\label{eq:t2}
\end{align}
\end{lemma}

\begin{proof}
This follows from an elementary calculation, which we briefly recall here.
For $d=1$, the claim follows by noting that $x_n = n/N$, and using the identity
\begin{gather}
\label{eq:trigid}
\sum_{\ell=0}^{N-1} q^\ell = 
\begin{cases}
\frac{q^N - 1}{q-1}, & (q \ne 1), \\
N, & (q = 1),
\end{cases}
\end{gather}
with $q = e^{2\pi i (m-n)/N}$ and $q = e^{2\pi i(k-k')/N}$, respectively. Indeed, assuming $d=1$ and denoting $-K := \min [[N]]$, then the above identity implies, for example,
\[
\sum_{k\in [[N]]} e^{2\pi i k(x_m-x_n)}
= \sum_{k \in [[N]]} \big[ \underbrace{e^{2\pi i (m-n)/N}}_{=:q} \big]^k
= \sum_{k \in [[N]]} q^k
= q^{-K} \sum_{\ell=0}^{N-1} q^\ell.
\]
If $q\ne 1$, then $q^N = e^{2\pi i (m-n)} = 1$. By \eqref{eq:trigid}, this implies that the last sum is $0$. On the other hand, if $q = 1$, then the last sum is trivially $=N$. We finally note that, for $m,n\in [N]$, we have $q=1$ if and only if $m=n$, implying that
\[
q^{-K} \sum_{\ell=0}^N q^\ell = N \delta_{mn}.
\]
Thus,
\[
\sum_{k\in [[N]]} e^{2\pi ik(x_m-x_n)} = N \delta_{mn},
\]
and \eqref{eq:t1} follows. The argument for \eqref{eq:t2} is analogous.
For $d>1$, the sum over $[[N]]^d = [[N]]\times \dots \times [[N]]$ is split into sums along each dimension, and the same argument is applied for each of the $d$ components, yielding the claim also for $d>1$.
\end{proof}

A \emph{trigonometric polynomial} $p: \mathbb{T}^d \mapsto \mathbb{R}^m$
is a function of the form
\begin{equation}
\label{eq:TP}
p(x) = \sum_{k\in [[N]]^d} c_k e^{2\pi i\langle k, x\rangle}
\end{equation}
with $c_k \in \mathbb{C}^m$ chosen to make $p(x)$ $\mathbb{R}^m$-valued at each 
$x \in \mathbb{T}^d.$ We note that the discrete and continuous $L^2$-norms 
are equivalent for trigonometric polynomials:

\begin{lemma}
\label{lem:polyequiv}
Let $N$ be a positive integer.
If $p(x)$ is a trigonometric polynomial, then 
\[
\frac1{N^{d/2}} \Vert p \Vert_{\ell^2(n\in [N]^d)}
=
\Vert p \Vert_{L^2(\Td)}.
\]
\end{lemma}

\begin{proof}
We have
\begin{align*}
\Vert p \Vert_{L^2(\Td)}^2
= \int_{\Td} |p(x)|^2 \, dx
= \sum_{k,k' \in [[N]]^d} c_k \overline{c}_{k'} \underbrace{\int_{\Td} e^{2\pi i\langle k-k', x\rangle} \, dx}_{=\delta_{kk'}}
= \sum_{k\in [[N]]^d} |c_k|^2,
\end{align*}
and
\begin{align*}
\frac{1}{N^d} \Vert p \Vert_{\ell^2(n\in [N]^d)}^2
&=
\frac1{N^d} \sum_{n\in [N]^d} |p(x_n)|^2
\\
&=
\sum_{k,k'\in [[N]]^d} c_k \overline{c}_{k'} \underbrace{\frac1{N^d} \sum_{n\in [N]^d} e^{2\pi i\langle k-k', x_n\rangle} }_{=\delta_{kk'}}
\\
&= \sum_{k\in [[N]]^d} |c_k|^2.
\end{align*}
This proves the claim.
\end{proof}

Let $v: \Td \to \R$ be a function with grid values $\{v(x_n)\}_{n\in [N]^d}$.
Let $\DFT(v)(k)$ denote the coefficients of the discrete Fourier transform
defined by \eqref{eq:DFTd}. Then 
\begin{align}
\label{eq:interp}
p(x) := \sum_{k\in [[N]]^d} \DFT(v)(k) e^{2\pi i\langle k, x\rangle},
\end{align}
is the trigonometric polynomial associated to $v$. 
The next lemma shows that $p(x)$ interpolates $v(x)$.

\begin{lemma}\label{lem:trig_grid_points}
The trigonometric polynomial $p(x)$ defined by \eqref{eq:interp} interpolates $v(x)$ at the grid points, i.e., we have $p(x_n) = v(x_n)$ for all $n\in [N]^d$.
\end{lemma}

\begin{proof}
Fix $n\in [N]^d$. Then 
\begin{align*}
p(x_n)
&= \sum_{k\in [[N]]^d} \DFT(v)(k) e^{2\pi i\langle k, x_n\rangle}
\\
&= \sum_{k\in [[N]]^d} \left\{\frac1{N^d} \sum_{m\in [N]^d} v(x_m) e^{-2\pi i \langle k, x_m\rangle} \right\} e^{2\pi i\langle k, x_n\rangle}
\\
&= \sum_{m\in [N]^d} v(x_m) \left\{
\frac1{N^d} \sum_{k\in [[N]]^d} e^{2\pi i \langle k, x_n - x_m \rangle }
\right\}
\\
&= \sum_{m\in [N]^d} v(x_m) \delta_{mn} 
\\
&= v(x_n),
\end{align*}
where we have made use of \eqref{eq:t1} to pass to the fourth line.
\end{proof}

The following trigonometric polynomial interpolation estimate 
for functions in Sobolev spaces $H^s(\Td)$ will be useful in
stating our main proposition:

\begin{lemma}
\label{lem:interpolation}
Let $v \in H^s(\Td)$ for $s>d/2$. Let $p$ denote the interpolating trigonometric polynomial given by \eqref{eq:interp}. Then 
\begin{align}
\label{eq:alias}
v(x) - p(x) 
= 
\sum_{k \in \Z^d \setminus [[N]]^d} \widehat{v}(k) e^{2\pi i \langle k, x \rangle } 
- \sum_{k \in [[N]]^d} \left\{
\sum_{\ell\in \Z^d\setminus \{0\}} \widehat{v}(k+\ell N)
\right\} e^{2\pi i\langle k, x\rangle }.
\end{align}
Furthermore, there exists a constant $c_{s,d}>0$, such that
\begin{align}
\label{eq:aliaserr}
\Vert v - p \Vert_{L^2(\Td)}
\le 
c_{s,d} \Vert v \Vert_{H^s(\Td)} N^{-s}.
\end{align}
\end{lemma}

\begin{remark}
The first sum on the right-hand side of \eqref{eq:alias} is the $L^2$-orthogonal Fourier projection of $v$ onto the complement of $\mathrm{span}\{ e^{2\pi i\langle k, x\rangle } \, |\, k\in [[N]]^d \}$. The second sum in \eqref{eq:alias} is known as an ``aliasing'' error; It arises because two Fourier modes are indistinguishable on the discrete grid whenever $k - k' \in N \Z^d$, i.e. $e^{2\pi i\langle k, x_n\rangle} = e^{2\pi i \langle k', x_n \rangle}$ for all $n\in [N]^d$. 
\end{remark}

\begin{proof}
Since $v\in H^s(\Td)$ has Sobolev smoothness $s$ for $s>d/2$, it can be shown that the Fourier series of $v$ is uniformly convergent. In particular, we may write 
\[ v(x) = \sum_{k' \in \Zd} \widehat{v}(k')e^{2\pi i \la k', x \ra}
\]
for 
\[
\widehat{v}(k') = \int_{\Td}v(x) e^{-2\pi i \la k', x \ra} \dx.
\]
First, substitution of $v(x_n) = \sum_{k'\in \Z^d} \widehat{v}(k') e^{2\pi i \langle k', x_n\rangle}$ into $\DFT(v)(k)$ yields
\begin{align*}
\DFT(v)(k)
&= \frac1{N^d} \sum_{n\in [N]^d} \left\{
\sum_{k'\in \Z^d} \widehat{v}(k') e^{2\pi i \langle k', x_n\rangle} 
\right\} e^{-2\pi i \langle k, x_n \rangle }
\\
&= \sum_{k'\in \Z^d} \widehat{v}(k')  \left\{ \frac1{N^d} \sum_{n\in [N]^d} 
e^{2\pi i \langle k'-k, x_n\rangle} 
\right\} 
\end{align*}
We now note that 
\[
\frac1{N^d} \sum_{n\in [N]^d} 
e^{2\pi i \langle k'-k, x_n\rangle} 
= 
\begin{cases}
0, &(k'\not \equiv k \mod{N}), \\
1, &(k' \equiv k \mod{N}),
\end{cases}
\]
as a consequence of the trigonometric identity \eqref{eq:t2}.
Letting $k' = k + \ell N$, i.e. $k'$ for which the sum inside the braces does not vanish, it follows that 
\[
\DFT(v)(k)
= \sum_{\ell \in \Z^d} \widehat{v}(k+\ell N).
\]
Thus,
\begin{align*}
v(x) - p(x) 
&= \sum_{k\in \Z^d} \widehat{v}(k) e^{2\pi i \langle k, x\rangle}
- \sum_{k\in [[N]]^d} \DFT(v)(k) e^{2\pi i \langle k, x\rangle}
\\
&= \sum_{k\in \Z^d\setminus [[N]]^d} \widehat{v}(k) e^{2\pi i \langle k, x\rangle}
+ \sum_{k\in [[N]]^d} \left\{ \widehat{v}(k) - \DFT(v)(k) \right\} e^{2\pi i \langle k, x\rangle}
\\
&= \sum_{k\in \Z^d\setminus [[N]]^d} \widehat{v}(k) e^{2\pi i \langle k, x\rangle}
- \sum_{k\in [[N]]^d} \left\{ \sum_{\ell \in \Z^d\setminus\{0\}} \widehat{v}(k+\ell N) \right\} e^{2\pi i \langle k, x\rangle}.
\end{align*}
We proceed to bound the last two terms. For the first term, we have by Parseval's theorem,
\begin{align*}
\Big\Vert 
\sum_{k\in \Z^d\setminus [[N]]^d} \widehat{v}(k) e^{2\pi i \langle k, x\rangle}
\Big\Vert_{L^2(\Td)}^2
&= \sum_{k\in \Z^d \setminus [[N]]^d} |\widehat{v}(k)|^2 
\\
&\le \frac{1}{(1+(N/2)^{2s})} \sum_{k\in \Z^d} (1+|k|^{2s}) |\widehat{v}(k)|^2 
\\
&\le 4^s N^{-2s} \Vert v \Vert^2_{H^s(\Td)},
\end{align*}
where $\Vert v \Vert^2_{H^s(\Td)} = \sum_{k\in \Z^d} (1+ |k|^{2s}) |\widehat{v}(k)|^2$, and for the second term
\begin{align*}
\Big\Vert 
\sum_{k\in [[N]]^d} \Big\{ \sum_{\ell \in \Z^d\setminus\{0\}} &\widehat{v}(k+\ell N) \Big\} e^{2\pi i \langle k, x\rangle}
\Big\Vert_{L^2(\Td)}^2
\\
&= \sum_{k\in [[N]]^d} \Big|\sum_{\ell \in \Z^d \setminus\{0\}} \widehat{v}(k+\ell N)\Big|^2 
\\
&\le
\sum_{k\in [[N]]^d} 
\Big(
\sum_{\ell \in \Z^d\setminus\{0\}} (1 + |k+\ell N|^{2s})^{-1}
\Big)
\\
&\qquad \times
\Big(
\sum_{\ell \in \Z^d\setminus\{0\}} (1 + |k+\ell N|^{2s}) |\widehat{v}(k+\ell N)|^2
\Big).
\end{align*}
The final inequality is obtained via Cauchy-Schwarz. We note that for $k\in [[N]]^d$, we have $|k|_\infty \le N/2$, and hence, for any integer vector $\ell \ne 0$, we obtain
\begin{equation}\label{eqn:KlN}
|k+\ell N|
\ge 
|k+\ell N|_\infty
\ge |\ell|_\infty N - |k|_\infty \ge |\ell|_\infty N - \frac{N}{2}
\ge \frac{N}{2} |\ell|_\infty \geq \frac{N}{2\sqrt{d}}|\ell|.
\end{equation}
We can now bound 
\begin{subequations}\label{eqn:CDS}
\begin{align}
\sum_{\ell \in \Z^d\setminus\{0\}} (1 + |k+\ell N|^{2s})^{-1}
&\le 
\sum_{\ell \in \Z^d\setminus\{0\}} \left( \frac{N}{2\sqrt{d}}\right)^{-2s} |\ell|^{-2s}
\\
&\le
c_{d,s} N^{-2s},
\end{align}
\end{subequations}
where $c_{d,s} := (4d)^s\sum_{\ell \in \Z^d\setminus\{0\}}  |\ell|^{-2s} < \infty$ is finite, since $s>d/2$ implies that the last series converges. Substitution of this bound in the estimate above implies,
\begin{align*}
\Bigg\Vert 
\sum_{k\in [[N]]^d} \Big\{ \sum_{\ell \in \Z^d\setminus\{0\}} &\widehat{v}(k+\ell N) \Big\} e^{2\pi i \langle k, x\rangle}
\Bigg\Vert_{L^2(\Td)}^2
\\
&\le
c_{d,s} N^{-2s}
\sum_{k\in [[N]]^d} 
\left(
\sum_{\ell \in \Z^d\setminus\{0\}} (1 + |k+\ell N|^{2s}) |\widehat{v}(k+\ell N)|^2
\right)
\\
&\le 
c_{d,s} N^{-2s} \Vert v \Vert_{H^s(\Td)}^2.
\end{align*}
Combining the above estimates, we conclude that 
\[
\Vert v - p \Vert_{L^2} 
\le c_{d,s} \Vert v \Vert_{H^s(\Td)} N^{-s},
\]
where we have re-defined $c_{d,s} := 2^s + (4d)^{s/2} \sum_{\ell\in \Z^d\setminus \{0\}} |\ell|^{-2s}$.
\end{proof}

We can now state the main outcome of this section:

\begin{proposition}
\label{prop:vp}
Let $v\in H^s(\Td)$ be given for $s>d/2$ and let  $\{v^N(x_n)\}_{n\in [N]^d}$ be any grid values. Let $p^N(x) = \sum_{k\in [[N]]^d} \DFT(v^N)(k) e^{2\pi i\langle k, x\rangle}$ be the interpolating trigonometric polynomial of $v^N$. Then,
\[
\Vert v - p^N \Vert_{L^2(\Td)}
\le \frac{1}{N^{d/2}} \Vert v- v^N \Vert_{\ell^2(n\in [N]^d)}
+ c_{d,s} \Vert v \Vert_{H^s(\Td)} N^{-s}.
\]
\end{proposition}

\begin{proof}
Let $p(x) = \sum_{k\in [[N]]^d} \DFT(v)(k) e^{2\pi i\langle k, x\rangle}$ be the interpolating trigonometric polynomial given the point-values $\{ v(x_n) \}_{n\in [N]^d}$. Then, 
\begin{align}
\label{eq:vp}
\Vert v - p^N \Vert_{L^2(\Td)}
\le 
\Vert v - p \Vert_{L^2(\Td)} + \Vert p - p^N \Vert_{L^2(\Td)}.
\end{align}
By Lemma \ref{lem:interpolation}, we have 
\[
\Vert v - p \Vert_{L^2(\Td)} \le c_{d,s} \Vert v \Vert_{H^s(\Td)} N^{-s}.
\]
By Lemma \ref{lem:polyequiv}, and since $p(x_n) = v(x_n)$, $p^N(x_n) = v^N(x_n)$ by Lemma \ref{lem:trig_grid_points}, we have
\begin{align*}
\Vert p - p^N \Vert_{L^2(\Td)}
&=
\frac{1}{N^{d/2}} \Vert p(x_n) - p^N(x_n) \Vert_{\ell^2(n\in [N]^d)}
\\
&=
\frac{1}{N^{d/2}} \Vert v(x_n) - v^N(x_n) \Vert_{\ell^2(n\in [N]^d)}.
\end{align*}
Substitution in \eqref{eq:vp} gives the claimed bound.
\end{proof}

\section{Discretization error derivation}\label{apdx:DiscErr}

In this section, we derive the error breakdown within each FNO layer. This error breakdown is used in the proofs of subsequent sections.  Within a single layer, we define the following quantities to track the error origin and propagation,
noting that, for values of $m_t$ that will vary with layer $t$, 
$\mcE_t^{(j)}: \XN \to \mathbb{R}^{m_t}$, $j=0,3$ and $\mcE_t^{(j)}: [[K]]^d \to \mathbb{C}^{m_t}$, $j=1,2$.

\begin{enumerate}[topsep=1.67ex,itemsep=0.5ex,partopsep=1ex,parsep=1ex]
\setcounter{enumi}{-1}
    \item $  \mcE^{(0)}_t(x_n) = v^N_{t}(x_n) - v_t(x_n) , \quad x_n \in 
\XN .$
    \item $ \mcE_t^{(1)}(k) = \frac{1}{N^d} \sum_{n\in [N]^d} v_t(x_n) e^{-2\pi i \la k, x_n \ra} - \int_{\Td} v_{t}(x) e^{-2\pi i\la k, x \ra} \; \dx, \quad 
k \in [[K]]^d.$
    \item $  \mcE^{(2)}_{t}(k) =  \frac{1}{N^d} \sum_{n\in[N]^d} \mcE^{(0)}_{t}(x_n)e^{-2\pi i \la k, x_n \ra}, \quad k \in [[K]]^d.$
    \item $    \mcE^{(3)}_{t}(x_n) = \sum_{k \in [[K]]^d} P_t^{(k)} \left(\mcE^{(1)}(k) + \mcE^{(2)}(k)\right)e^{2\pi i \langle k, x_n \rangle}, \quad x_n \in \XN.$
    \item $    \mcE^{(0)}_{t+1}(x_n) =\sigma\left(W_tv_t(x_n) + \mcK_tv_t(x_n) + b_t + W_t \mcE^{(0)}_t(x_n)+ \mcE^{(3)}_t(x_n)\right)\\ - \sigma(W_tv_t(x_n) + \mcK_t v_t(x_n) + b_t)= v_{t+1}^N(x_n) -v_{t+1}(x_n), \quad x_n \in \XN.$
\end{enumerate}
Here, $\mcEzero_t$ is the initial error in the inputs to FNO layer $t$, $\mcEone$ is the aliasing error, $\mcEtwo_t$ is the initial error $\mcEzero_t$ after the discrete Fourier transform, and $\mcEthree_t$ is the error after the operation of the kernel $\mcK_t$. Finally, the initial error for the next layer is given by $\mcEzero_{t+1}$ in terms of the error quantities of the previous layer. Intuitively, the quantity $\mcEone$ is the source of the error within each layer since it depends only on the ground truth $v_t$. All other error quantities are propagation of existing error from previous layers. We provide an exact derivation of these quantities in the following. 

Let $\mcE^{(0)}_{t}$ be the error in the inputs to FNO layer $t$ such that \begin{equation*}
    \mcE^{(0)}_t(x_n) = v^N_{t}(x_n) - v_t(x_n) , \quad x_n \in \XN. 
\end{equation*}
Let $\mcF(v_t)(k) = \int_{\Td} v_t(x) e^{-2\pi i \la k, x \ra} \dx$ denote the Fourier transform and $\DFT$ as in equation \eqref{eqn:DFT}. Then for $k \in [[K]]^d,$
\begin{align*}
    \DFT(v^N_{t})(k) &= \frac{1}{N^d} \sum_{n\in [N]^d} v_{t}(x_n)e^{-2\pi i\la k, x_n \ra} + \frac{1}{N^d} \sum_{n\in [N]^d} \mcE^{(0)}_{t}(x_n)e^{-2\pi i \la k, x_n \ra} \\
    &= \mcF(v_{t})(k)  + \mcE^{(1)}_{t}(k) + \mcE^{(2)}_{t}(k)
\end{align*}
where $\mcE^{(1)}_t$ is the error resulting from computing the Fourier transform of $v_t$ on a discrete grid rather than all of $\Td$, i.e.
\begin{equation*}
    \mcE_t^{(1)}(k) = \frac{1}{N^d} \sum_{n\in [N]^d} v_t(x_n) e^{-2\pi i \la k, x_n \ra} - \int_{\Td} v_{t}(x) e^{-2\pi i\la k, x \ra} \; \dx
\end{equation*}
and $\mcE^{(2)}_t$ is the error $\mcE^{(0)}_t$ after the discrete Fourier transform, i.e.
\begin{equation*}
    \mcE^{(2)}_{t}(k) =  \frac{1}{N^d} \sum_{n\in[N]^d} \mcE^{(0)}_{t}(x_n)e^{-2\pi i \la k, x_n \ra}.
\end{equation*}
For $x_n \in \XN$, the output of the discrete kernel integral operator acting on $v_t^N$ is given by 
\begin{align*}
    (\mcK^N_tv^N_t)(x_n) &= \sum_{k \in [[K]]^d}P_t^{(k)} \left(\mcF(v_{t})(k) + \mcE^{(1)}_{t}(k) + \mcE^{(2)}_{t}(k)\right)e^{2\pi i \langle k, x_n \rangle}\\
    &= (\mcK_t v_t)(x_n) + \mcE^{(3)}_{t}(x_n)
\end{align*}
where 
\begin{equation*}
    \mcE^{(3)}_{t}(x_n) = \sum_{k \in [[K]]^d} P_t^{(k)} \left(\mcE^{(1)}(k) + \mcE^{(2)}(k)\right)e^{2\pi i \langle k, x_n \rangle}.
\end{equation*}

Finally, the output of layer $t$ is given by 
\begin{align*}
    v^N_{t+1}(x_n) & = \sigma\left(W_t\bigl(v_t(x_n) + 
\mcE^{(0)}_{t}(x_n)\bigr) + (\mcK^N_t v^N_t)(x_n) + b_t\right)\\
    & = \sigma\left(W_tv_t(x_n) + \mcK_tv_t(x_n) + b_t + W_t \mcE^{(0)}_t(x_n)+ \mcE^{(3)}_t(x_n)\right). 
\end{align*}
Therefore, the initial error for the next layer is given by 
\begin{align*}
    \mcE^{(0)}_{t+1}(x_n) &=\sigma\left(W_tv_t(x_n) + \mcK_tv_t(x_n) + b_t + W_t \mcE^{(0)}_t(x_n)+ \mcE^{(3)}_t(x_n)\right) \\ &- \sigma\left(W_tv_t(x_n) + \mcK_t v_t(x_n) + b_t\right).
\end{align*}

\section{Proofs of approximation theory lemmas} \label{appx:lemmas}
 We bound the components described in Appendix \ref{apdx:DiscErr} in the following proposition.

\begin{proposition} \label{prop:bndstep} Under Assumptions \ref{asst:main}, it holds that 
\begin{enumerate}
    \item $\|\mcEone_t\|_{\ell^2(k \in [[K]]^d)} \leq \alpha_{d,s}N^{-s}\|v_t\|_{H^s}$ where $\alpha_{d,s}$ is independent of $N,v_t$;
    \item $\|\mcEtwo_t\|_{\ell^2(k \in [[N]]^d)} = N^{-d/2}\|\mcEzero_t\|_{\ltwoN}$;
    \item $\|\mcEthree_t\|_{\ell^2(n \in [N]^d)} \leq N^{d/2}\|P_t\|_F\left(\|\mcEone_t\|_{\ell^2(k \in [[K]]^d)} + \|\mcEtwo_t\|_{\ell^2(k\in [[K]]^d)}\right);$
    \item $\|\mcEzero_{t+1}\|_{\ltwoN} \leq B\left(\|W_t\|_2 \|\mcEzero_t\|_{\ltwoN} + \|\mcEthree_t\|_{\ltwoN}\right)$.
\end{enumerate}
\end{proposition}
\begin{proof}
    Beginning with the definition of $\mcEone_t(k)$, we have
    \begin{align*}
        \|\mcEone_t\|_{\ell^2(k \in [[K]]^d)}^2 & = \Big\|\frac{1}{N^d} \sum_{n \in [N]^d} v_t(x_n) e^{-2\pi i \la k, x_n \ra} - \int_{\Td}e^{-2\pi i \la k, x \ra}v_t(x)\dx \Big\|_{\ell^2(k \in [[K]]^d)}^2.
    \end{align*}
    Denote the terms in the above expression $\vhat_t^N(k)$ and $\vhat_t(k)$, respectively. Since $s > \frac{d}{2}$,
    \begin{equation*}
        v_t(x_n) = \sum_{k \in \Zd}\vhat_t(k)e^{2\pi i\la k, x_n\ra},
    \end{equation*}and it follows that 
    \begin{align*}
        \vhat^N_t(k') & = \frac{1}{N^d} \sum_{n \in [N]^d}\left(\sum_{k \in \Zd}\vhat_t(k) e^{2\pi i \la k, x_n \ra}\right)e^{-2\pi i \la k', x_n \ra}\\
        & = \sum_{k \in \Zd}\vhat_t(k) \frac{1}{N^d}\sum_{n \in [N]^d}e^{2\pi i \la k - k', x_n \ra}\\
        & = \sum_{\ell \in \Zd}\vhat_t(k' + \ell N) .
    \end{align*}
    Therefore, 
    \begin{align*}
    \|\mcEone_t&\|_{\ell^2(k \in [[K]]^d)}^2  = \|\vhat^N_t- \vhat_t\|^2_{\ell^2(k \in [[K]]^d)}\\
    & = \sum_{k \in [[K]]^d}\left|\sum_{\ell \in \Zd\setminus\{0\}} \vhat_t(k + \ell N)\right|^2\\
    & \leq \sum_{k \in [[K]]^d}\left(\sum_{\ell \in \Zd\setminus\{0\}}\frac{1}{|k + \ell N|^{2s}}\right)\sum_{\ell \in \Zd \setminus\{0\}}|k + \ell N|^{2s} |\vhat_t(k + \ell N)|^2
    \end{align*}
    by Cauchy-Schwarz. We bound each component separately. It is clear from 
Definition \ref{eqn:Hs} that 
    \begin{equation}
        \sum_{k \in [[K]]^d} \sum_{\ell \in \Zd \setminus\{0\}}|k + \ell N|^{2s} |\vhat_t(k + \ell N)|^2\leq \|v_t\|_{H^s}^2.
    \end{equation}
    To bound the first component independently of $k$, we note from $K \leq \frac{N}{2}$ and equation \eqref{eqn:KlN} that 
    \begin{align*}
        \sum_{\ell \in \Z^d\setminus \{0\}}\frac{1}{|k + \ell N|^{2s}}& \leq \sum_{\ell \in \Z^d\setminus \{0\}}\left(\frac{N}{2\sqrt{d}} |\ell|\right)^{-2s} \\
        & \leq \alpha_{d,s}^2 N^{-2s}
    \end{align*}
    by equation \eqref{eqn:CDS}, where $\alpha^2_{d,s} = (4d)^s\sum_{\ell \in \Z^d\setminus \{0\}} |\ell|^{-2s}$ is finite since $s \geq \frac{d}{2}$. We express the final bound as 
    \begin{equation*}
        \|\mcEone_t\|_{k \in [[K]]^d} \leq \alpha_{d,s}N^{-s}\|v_t\|_{H^s}.
    \end{equation*}
    For $\mcEtwo_t(k)$ we have the definition 
    \begin{equation*}
        \mcEtwo_t(k) = \frac{1}{N^d}\sum_{n \in [N]^d} \mcEzero_t(x_n)e^{-2\pi i \la k, x_n \ra}.
    \end{equation*}
    By Parseval's Theorem, we have 
    \begin{equation}
        \|\mcEtwo_t\|^2_{\ell^2(k \in [[N]]^d)} = \frac{1}{N^d} \|\mcEzero_t\|^2_{\ltwoN}.
    \end{equation}
    For $P_t \in \mathbb{R}^{d_{v_{t+1}}\times K^d \times d_{v_t}}$ we define the tensor Frobenius norm 
    \newline $\|P_t\|^2_F = \sum_{k \in [[K]]^d} \|P_t^{(k)}\|_F^2$.
    \begin{align*}
        \|\mcEthree_t\|_{\ell^2(n \in [N]^d)}^2 & = \sum_{n \in [N]^d} \left | \sum_{k \in [[K]]^d}\ P_t^{(k)} \left(\mcEone_t(k) + \mcEtwo_t(k)\right)e^{2\pi i \la k, x_n \ra}\right|^2\\
        & \leq N^d \left|\sum_{k \in [[K]]^2} |P_t^{(k)}(\mcEone_t(k) + \mcEtwo_t(k))|\right|^2\\
        & \leq N^d \sum_{k \in [[K]]^d} \|P_t^{(k)}\|_F^2 \sum_{k \in [[K]]^d}|\mcEone_t(k) + \mcEtwo_t(k)|^2\\
        &= N^d\|P_t\|_F^2\|\mcEone_t+ \mcEtwo_t\|^2_{\ell^2(k \in [[K]]^d)}\\
        \|\mcEthree_t\|_{\ell^2(n \in [N]^d)}  &\leq N^{d/2}\|P_t\|_F\left(\|\mcEone_t\|_{\ell^2(k \in [[K]]^d)} + \|\mcEtwo_t\|_{\ell^2(k\in [[K]]^d)}\right)
    \end{align*}
    Finally, we have
    \begin{align*}
        \|\mcEzero_{t+1}\|_{\ltwoN}^2 & = \sum_{n \in [N]^d}\left|\sigma(W_tv_t + \mcK_tv_t + b_t + W_t\mcEzero_t(x_n) + \mcEthree_t(x_n)) - \sigma(W_tv_t + \mcK_tv_t + b_t)\right|^2\\
        & \leq \sum_{n\in[N]^d}B^2\left|W_t\mcEzero_t(x_n) + \mcEthree_t(x_n)\right|^2\\
        \|\mcEzero_{t+1}\|_{\ltwoN} & \leq B\left(\|W_t\|_2 \|\mcEzero_t\|_{\ltwoN} + \|\mcEthree_t\|_{\ltwoN}\right)
    \end{align*}
    where $\|\cdot\|_2$ is the matrix-$2$ norm. Recall $B$ bounds derivatives of $\sigma$ in Assumptions \ref{asst:main}.

\end{proof}
The results of Proposition \ref{prop:bndstep} allow us to easily prove the following lemma.
\begin{restatable}{lemmma}{lemOneLayer}
        \label{lem:layerresult} Under Assumptions \ref{asst:main}, the following bound holds: 
\begin{equation}
\frac{1}{N^{d/2}}\|\mcEzero_{t+1}\|_{\ltwoN} \leq BM \left(\frac{2}{N^{d/2}} \|\mcEzero_t\|_{\ltwoN} + \alpha_{d,s}N^{-s}\|v_t\|_{H^s}\right)
\end{equation}
where $\alpha_{d,s}$ is a constant dependent only on $d$ and $s$.
\end{restatable}
\begin{proof}
    From Proposition \ref{prop:bndstep}, and shortening the notation $\ltwoN$ to $\ell^2$,  
    \begin{align*}
        \|\mcEzero_{t+1}\|_{\ell^2} & \leq B\left(\|W_t\|_2\|\mcEzero_t\|_{\ell^2}+ N^{d/2}\|P_t\|_F\left(\alpha_{d,s}N^{-s}\|v_t\|_{H^s} + N^{-d/2}\|\mcEzero_t\|_{\ell^2}\right)\right)
    \end{align*}
    Combining terms gives
    \begin{equation}
\|\mcEzero_{t+1}\|_{\ell^2}  \leq B\left(\big(\|W_t\|_2 + \|P_t\|_F\big) \|\mcEzero_t\|_{\ell^2} + \alpha_{d,s}N^{d/2-s} \|P_t\|_F\|v_t\|_{H^s}\right).
\end{equation}
Replacing $\|W_t\|_2$ and $\|P_t\|_F$ with $M$ and rescaling gives 
\begin{equation*}
    \frac{1}{N^{d/2}}\|\mcEzero_{t+1}\|_{\ltwoN} \leq BM \left(\frac{2}{N^{d/2}} \|\mcEzero_t\|_{\ltwoN} + \alpha_{d,s}N^{-s}\|v_t\|_{H^s}\right).
\end{equation*}
\end{proof}

\section{Proofs of regularity theory lemmas} \label{appx:lemmas2}

The proof of Lemma \ref{lem:HsBound} relies on another result for bounding the $H^s$ norm of compositions of functions, which is largely taken from the lemma of \citet[sec. 2~, p.~273]{moser1966rapidly} without assuming an $L^{\infty}$ norm of $v$ less than $1$. We state a proof here for completeness.

\begin{lemma}\label{lem:Moser}
Assume $\varphi: \Td \to \Td$ possesses continuous derivatives up to order $r$ which are bounded by $B$. Then 
\begin{equation*}
    |\varphi \circ v|_r \leq Bc \left(1+ \|v\|_{\infty}^{r-1}\right)\|v\|_{H^r}
\end{equation*}
provided $v \in H^r(\Td)$, where $c$ is a constant dependent on $r$ and $d$. 
\end{lemma}
\begin{proof}
    By Fa\`a di Bruno's formula, we have 
    \begin{equation}\label{eqn:FdB}
        D_x^r(\varphi \circ v(x))  = \sum C_{\alpha, r}\frac{d^{\rho}\varphi}{d x^{\rho}}(v(x)) \prod_{j=1}^r (D_x^j v(x))^{\alpha_j}
    \end{equation}
    where the sum is over all nonnegative integers $\alpha_1, \dots, \alpha_r$ such that $\alpha_1 + 2\alpha_2 + \dots + r\alpha_r = r$, the constant $C_{\alpha, r} = \frac{r!}{\alpha_1! \alpha_2! 2!^{\alpha_2}\dots \alpha_r!r!^{\alpha_r}}$, and $\rho : = \alpha_1 + \alpha_2 + \dots + \alpha_r$. 

    We seek a bound on square integrals of \eqref{eqn:FdB}. Setting $v_0 = \frac{d^{\rho}\varphi}{d x^{\rho}}v$, $v_{\lambda} = D_x^{\lambda}v$, $\alpha_0 = 1$, $p_0 = \infty$, and $p_{\lambda} = \frac{r}{\lambda \alpha_{\lambda}}$ and noting that $\sum_{\lambda = 0}^r \frac{1}{2p_{\lambda}} = \frac{1}{2}$,  we have by H\"older's inequality for multiple products that
    \begin{align*}
        \int_{\Td} \left|\frac{d^{\rho}\varphi}{d x^{\rho}}(v(x)) \prod_{j=1}^r (D_x^j v(x))^{\alpha_j}\right|^2 \dx \leq \int_{\Td}\prod_{\lambda = 0}^r |v_{\lambda}|^{2\alpha_{\lambda}} \dx &\leq \prod_{\lambda = 0}^r \left(\int_{\Td} |v_{\lambda}|^{2\alpha_{\lambda}p_{\lambda}} \dx \right)^{1/{p_{\lambda}}}\\
        & = \|v_0\|_{\infty}^2 \prod_{\lambda = 1}^r \left(\int_{\Td} |v_{\lambda}|^{2\alpha_{\lambda}p_{\lambda}}\dx\right)^{1/{p_{\lambda}}}
    \end{align*}
    The first factor is bounded above by $B^2$ by assumption. By application of Gagliardo-Nirenberg, the second factor may be bounded by 
    \begin{align*}
        \prod_{\lambda = 1}^r \left(\int_{\Td}|D_x^{\lambda}v|^{2r/\lambda}\dx\right)^{\lambda \alpha_{\lambda}/r} &\leq C^r \prod_{\lambda = 1}^r \|v\|_{\infty}^{2\alpha_{\lambda}(1-\lambda/r)}\left(\|D_x^rv\|^2 + \|v\|^2\right)^{\alpha_{\lambda}\lambda/r}\\
        & \leq C^r \|v\|_{\infty}^{2\rho -2}\|v\|_{H^r}^2
    \end{align*}
      since $\sum_{\lambda} \lambda \alpha_{\lambda} = r$, and $\sum_{\lambda}\alpha_{\lambda} = \rho$. Combining the bounds, 
    \begin{equation*}
        \int_{\Td}\prod_{\lambda = 0}^r |v_{\lambda}|^{2\alpha_{\lambda}}\dx \leq B^2 C^r \|v\|_{\infty}^{2\rho -2}\|v\|_{H^r}^2.
    \end{equation*}
    If $\|v\|_{\infty} <1$, we have the bound 
    \begin{equation}
        \int_{\Td}\prod_{\lambda = 0}^r |v_{\lambda}|^{2\alpha_{\lambda}}\dx \leq B^2 C^r\|v\|_{H^r}^2,
    \end{equation}
    and otherwise since $\rho \leq r$,
    \begin{equation}\label{eqn:boundnot1}
        \int_{\Td}\prod_{\lambda = 0}^r |v_{\lambda}|^{2\alpha_{\lambda}}\dx \leq B^2 C^r\|v\|_{\infty}^{2r-2}\|v\|^2_{H^r}.
    \end{equation}
     Since these bounds hold for any term in the sum \ref{eqn:FdB}, we obtain
    \begin{equation}
        |\varphi \circ v|_r \leq Bc \left(1+ \|v\|_{\infty}^{r-1}\right)\|v\|_{H^r}
    \end{equation}
    for a different constant $c$ depending on $r$ and $d$. 
\end{proof}

Now we may prove Lemma \ref{lem:HsBound}.
\lemHs*
\begin{proof}
    First we bound $\|\mcK_tv_t\|_{\infty}$. Recall $\widehat{v}_t(k): = \int_{\Td} v_t(x) e^{-2\pi i \la k, x\ra} \dx$.
    \begin{align*}
        \|\mcK_tv_t\|_{\infty} & = \|\sum_{k \in [[K]]^d}P^{(k)}_t \widehat{v}_t(k)e^{2\pi i \la k, x\ra}\|_{\infty} \\
        & \leq \sum_{k \in [[K]]^d}\|P^{(k)}_t\||\widehat{v}_t(k)|\\
        & \leq \left(\sum_{k \in [[K]]^d} \|P_t^{(k)}\|^2\right)^{1/2}\|\widehat{v}_t\|_{\ell^2(k \in [[K]]^d)}\\
        & \leq \|P_t\|_FK^{d/2}\|\widehat{v}_t\|_{\ell^2(k \in [[K]]^d)}\\
        & \leq \|P_t\|_FK^{d/2} \|v_t\|_{L^2(\Td)}.
    \end{align*}
    Then 
    \begin{align*}
        \|W_tv_t + \mcK_tv_t + b_t\|_{\infty} \leq \|W_t\|_2\|v_t\|_{\infty} + |b_t| + \|P_t\|_FK^{d/2} \|v_t\|_{L^2(\Td)},
    \end{align*}
    and by Lipschitzness of $\sigma$ we have
    \begin{align*}
        \|v_{t+1}\|_{\infty} \leq \sigma^* + BM\left(1 + \|v_t\|_{\infty} + K^{d/2}\|v_t\|_{L^2(\Td)}\right).
    \end{align*}
    Next we bound $|v_{t+1}|_s$. Letting $f_t = W_tv_t + \mcK_tv_t + b_t$, we see from Lemma \ref{lem:Moser} that bounding $\|f_t\|_{H^s}$ will give the result. 
    \begin{align*}
        D_x^s(f_t) &= W_t(D_x^sv_t) + \mcK_t(D^s_x v_t).\\
        \int_{\Td}|D^s_x (f_t)|^2 \dx & \leq 2\left(\int_{\Td} |W_t (D^s_x v_t)|^2 \dx + \int_{\Td} |\mcK_t(D^s_x v_t)|^2 \dx \right)
    \end{align*}
    The first integral on the right may be bounded by $\|W_t\|_2^2|v_t|_s^2$. To bound the second integral, 
    \begin{align*}
        \int_{\Td}|\mcK_t(D^s_x v_t)|^2 \dx & = \int_{\Td} \left|\sum_{k \in [[K]]^d} P_t^{(k)} \widehat{g}_t(k) e^{2\pi i \la k, x\ra } \right|^2\dx
    \end{align*}
    where $\widehat{g}_t(k)$ are the Fourier coefficients of $D^s_x v_t$. Continuing, 
    \begin{align*}
         \int_{\Td}|\mcK_t(D^s_x v_t)|^2 \dx & \leq \int_{\Td}\|P_t\|_F^2 \sum_{k \in [[K]]^d}|\widehat{g}_t(k)|^2 \dx\\
         & \leq \|P_t\|_F^2\|D^s_xv_t\|_{L^2}^2,
    \end{align*}
    giving a bound of 
    \begin{equation*}
        |f_t|_s \leq 2M|v_t|_s
    \end{equation*}
    In the following, $\lesssim$ denotes inequality up to a constant multiple that does not depend on any of the variables involved. Combining Lemma \ref{lem:Moser} and the above bounds, we have 
    \begin{align*}
        |\sigma \circ f_t|_s & \leq Bc(1+\|f_t\|_{\infty}^{s-1})\|f_t\|_{H^s}\\
        & \leq Bc(1+ (M(1+\|v_t\|_{\infty} + K^{d/2}\|v_t\|_{\infty}))^{s-1})(M(1+\|v_t\|_{\infty}+K^{d/2}\|v_t\|_{\infty})+2M|v_t|_s)\\
        & \lesssim BcM^sK^{ds/2}(1+(1+\|v_t\|_{\infty})^{s-1})(1+\|v_t\|_{\infty}+|v_t|_s)\\
        & \lesssim BcM^sK^{ds/2}(1+\|v_t\|_{\infty})^{s-1}(1+\|v_t\|_{\infty})(1+|v_t|_s)\\
        & \lesssim BcM^sK^{ds/2}(1+\|v_t\|_{\infty})^s(1+|v_t|_s).
    \end{align*}
\end{proof}

\section{Proof of Theorem \ref{thm-fno:main}}

\label{appx:main}

\thmmain*
\begin{proof}
    Temporarily dropping the notational dependence of $v_t$ and $v_0$ on $a$, from Lemma \ref{lem:HsBound} we have for $t \geq 1$, 
    \begin{align*}
        \|v_t\|_{\infty} &\lesssim \sigma^*\sum_{j=0}^{t-1}(BMK^{d/2})^j + \sum_{j=1}^t (BMK^{d/2})^j+(BMK^{d/2})^t\|v_0\|_{\infty}\\
        |v_t|_s & \lesssim \left(\sum_{j = 1}^t (BcM^sK^{ds/2})^j \prod_{\ell = t-j}^{t-1}(1+\|v_{\ell}\|_{\infty})^s\right) + (BcM^sK^{ds/2})^t\left(\prod_{\ell = 0}^{t-1}(1+\|v_{\ell}\|_{\infty})^s\right)|v_0|_s. 
    \end{align*}
    Denote $\max\{BMK^{d/2},B^{1/s}c^{1/s}MK^{d/2},1\}$ by $C_0$. Since $\sigma^* \geq 1$, the bound on $\|v_t\|_{\infty}$ simplifies to 
    \begin{align*}
        \|v_t\|_{\infty} &\lesssim \sigma^* \sum_{j=1}^{t}C_0^j + C_0^t\|v_0\|_{\infty}\\
        &\leq \sigma^* tC_0^t + C_0^t \|v_0\|_{\infty}
    \end{align*}
    Plugging in this bound to the product in the bound on $|v_t|_s$, we have
    \begin{align*}
        \prod_{\ell = t-j}^{t-1} (1+\|v_{\ell}\|_{\infty})^s & \lesssim \prod_{\ell = t-j}^{t-1} (1+\ell\sigma^*C_0^{\ell} + C_0^{\ell}\|v_0\|_{\infty})^s\\
        & \lesssim C_0^{tsj}(t)^{sj}(\sigma^*+\|v_0\|_{\infty})^{sj}.
    \end{align*}
    Combining these two bounds, we attain the following bound on $|v_t|_{s}$ for $t \geq 1$. 
    \begin{align*}
        |v_t|_s & \lesssim \left(\sum_{j = 1}^t (C_0)^{sj} C_0^{tsj}(t)^{sj}(\sigma^*+\|v_0\|_{\infty})^{sj}\right) + C_0^{ts}\left(C_0^{t^2s}(t)^{st}(\sigma^*+\|v_0\|_{\infty})^{st}\right)|v_0|_s\\
        & \lesssim \left(\sum_{j=1}^t C_0^{2tsj}(t)^{sj}(\sigma^*+\|v_0\|_{\infty})^{sj} \right) + C_0^{2t^2s}(t)^{st}(\sigma^*+\|v_0\|_{\infty})^{st}|v_0|_s\\
        & \lesssim (C_0^{2t^2s}t^{st+1} +C_0^{2t^2s}t^{st}|v_0|_s)(\sigma^* + \|v_0\|_{\infty})^{st}
    \end{align*}

    and the following bound on $\|v_t\|_{H^s}$
    \begin{align}
    \label{eqn:Hsbound}
        \|v_t\|_{H^s} \lesssim (C_0^{2t^2s}t^{st+1}|v_0|_s)(\sigma^* + \|v_0\|_{\infty})^{st} + \sigma^* tC_0^t + C_0^t \|v_0\|_{\infty}.
    \end{align}
     \mt{Recall that $v_0 = \cP(a)$, and $\cP$ is a shallow neural network, which is a special case of a Fourier layer where the coefficients $P_t^{(k)}$ are set to $0$. Assumptions \ref{asst:main} include boundedness of the coefficients of $\cP$ by $M$. Thus we may increment $t$ by $1$ in the bound and write
     \begin{subequations}\label{eqn:Hsbound2}
    \begin{align}
        &\sup_{a\in \cA_c}\|v_t(a)\|_{H^s} \\&\lesssim \sup_{a \in \cA_c} (C_0^{2(t+1)^2s}(t+1)^{s(t+1)+1}|a|_s)(\sigma^* + \|a\|_{\infty})^{s(t+1)} + \sigma^* (t+1)C_0^{t+1} + C_0^{t+1} \|a\|_{\infty}.
    \end{align}
    \end{subequations}
    Since $\cA$ is a compact set in $H^s$, and $s > \frac{d}{2}$, both $\|a\|_{\infty}$ and $|a|_s$ are bounded uniformly over $\cA$ by a constant depending on $\cA$ since $H^s$ is continuously embedded in $L^{\infty}$. Thus, we may }denote this upper bound by $C_1$, which does not depend on $N$. Let $\mathcal{E}_{t+1}^{(0)}(a) = v^N_t(a) - v_t(a)$. Then from Lemma \ref{lem:layerresult}, we have 
    \begin{align*}
      \sup_{a \in \cA_c} \frac{1}{N^{d/2}}\|\mcEzero_{t+1}(a)\|_{\ltwoN} \lesssim BM \left(\frac{2}{N^{d/2}} \sup_{a \in \cA_c}\|\mcEzero_t(a)\|_{\ltwoN} + \alpha_{d,s}N^{-s}C_{1}\right).
    \end{align*}
    By the discrete Gronwall lemma, 
    \begin{align*}
        &\sup_{a\in\cA_c}\frac{1}{N^{d/2}}\|\mcEzero_t(a)\|_{\ltwoN} \\ &\lesssim \frac{BM\alpha_{d,s}N^{-s}C_{1}}{1-2BM}(1-(2BM)^t) +\frac{1}{N^{d/2}}\sup_{a\in \cA_c}\|\mcEzero_0(a)\|_{\ltwoN}(2BM)^t.
    \end{align*}
    Since we assume we begin with no error, $\|\mcEzero_0(a)\|_{\ltwoN} = 0$, this simplifies to 
    \begin{align*}
        \sup_{a \in \cA_c}\frac{1}{N^{d/2}}\|\mcEzero_t(a)\|_{\ltwoN} & \lesssim \frac{BM\alpha_{d,s}C_{1}}{1-2BM}(1-(2BM)^t)N^{-s}.
    \end{align*}
    Denoting the factor in front of $N^{-s}$ by $C$ and absorbing the effects of $\lesssim$ into $C$, we have the result that 
    \begin{equation*}
        \sup_{a \in \cA_c}\frac{1}{N^{d/2}}\|v_t(a) - v_t^N(a)\|_{\ell^2(n \in [N]^d)} \leq CN^{-s}.
    \end{equation*}
\end{proof}

\begin{remark}
    A trivial consequence of the above theorem is that under Assumptions \ref{asst:main},
\begin{equation*}\label{eqn:convergence}
     \lim_{N \to \infty} \sup_{a \in \cA_c} \frac{1}{N^{d/2}}\|v_t^N(a) - v_t(a)\|_{\ltwoN} = 0.
\end{equation*}
 Indeed, a stronger result holds  that the discrete $\ell^{\infty}$ norm converges at a rate $N^{-s+d/2}$ by a straightforward inverse inequality. 
\end{remark}

\section{Proof of Theorem \ref{thm:main'}}
\label{appx:main'}
\thmmainp*
\begin{proof}
We temporarily drop the dependence of $p^N_t$ and $v^N_t$ on $a$. Let $p^N_t(x)$ be the interpolating trigonometric polynomial associated with the data $\{v^N_t(x_n)\}_{n\in [N]^d}$. By Proposition \ref{prop:vp}, we have 
\[
\Vert v_t - p^N_t \Vert_{L^2(\Td)}
\le 
\frac{1}{N^{d/2}} \Vert v_t - v^N_t\Vert_{\ell^2(n\in [N]^d)}
+
c_{d,s} \Vert v \Vert_{H^s(\Td)} N^{-s}.
\]
By \eqref{eqn:Hsbound2}, we have $\sup_{a\in\cA_c}\Vert v_t(a) \Vert_{H^s(\Td)} \le C_{1}$.
Furthermore, it follows from Theorem \ref{thm-fno:main}, that 
\[
\sup_{a \in \cA_c} \frac{1}{N^{d/2}} \Vert v_t(a) - v^N_t(a) \Vert_{\ell^2(n\in [N]^d)}
\le 
C N^{-s}.
\]
We conclude that
\begin{align*}
\sup_{a\in\cA_c}\Vert v_t(a) - p^N_t(a) \Vert_{L^2(\Td)}
&\le 
(C + c_{d,s}C_1) N^{-s}
\end{align*}
Thus, the claimed bound holds with $C' = C + c_{d,s}C_1$.

\end{proof}

\section{Additional numerical results}\label{appd:add_num_res}

\begin{figure}[ht!]
    \begin{center}
    \includegraphics[width = \linewidth]{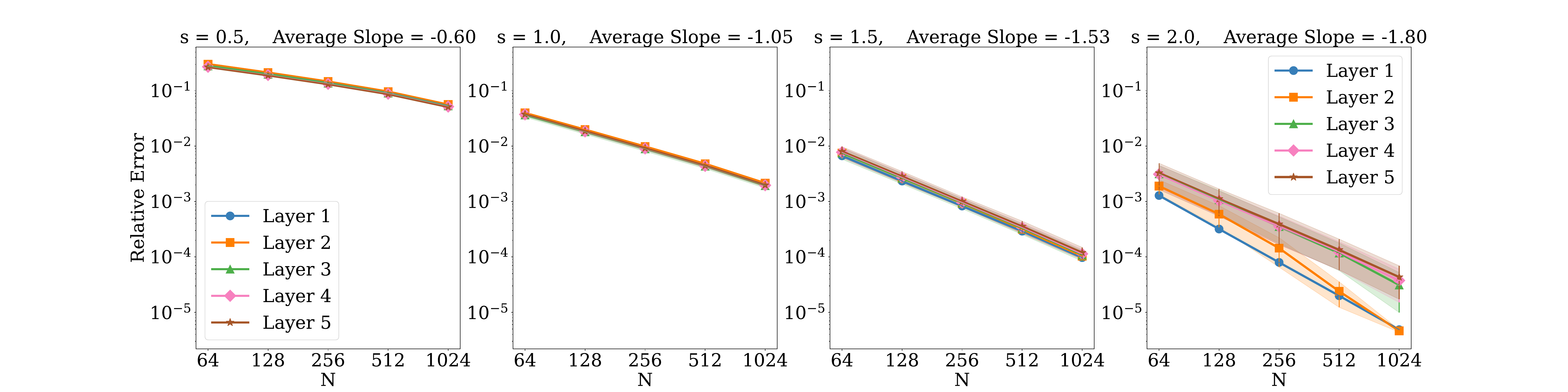}
    \caption{Relative error versus $N$ and $s$ for an FNO with default $\times 10$ initial weights.}
    \label{fig:err_vs_N_10x}
    \end{center}
\end{figure}

Figure \ref{fig:err_vs_N_10x} addresses the question of error decreasing or increasing with layer count. The figure shows that when the FNO weights are randomly initialized with the default initialization and then multiplied by $10$, the error increases with the number of layers instead of decreases. Additionally, in this model, the large weights mean that the GeLU activation acts like a ReLU activation for smaller discretizations. This phenomenon is apparent for inputs with regularity $s=2$, where the first layer has the appropriate slope, but the other layers only begin to approach that rate at higher discretizations. Earlier layers achieve this rate first because of the smaller magnitude state norm in earlier layers for this model.

\begin{figure}[ht!]
    \begin{center}
    \includegraphics[width = \textwidth]{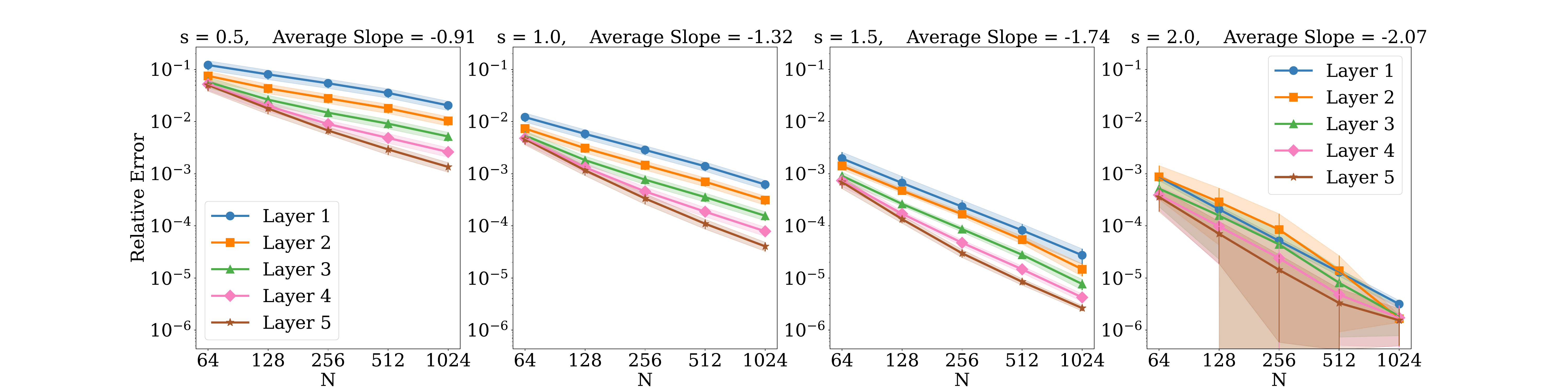}
    \caption{Relative error versus $N$ and $s$ for an FNO with all weights equals to $1$.}
    \label{fig:err_vs_N_all1s}
    \end{center}
\end{figure}

As an alternative setting of the weights, Figure \ref{fig:err_vs_N_all1s} shows the discretization error when all the weights are set to $1$. In this case, the error is more erratic. The error decreases faster than expected and with less consistency than the Gaussian weight models, and the decay rate increases with each layer. In this sense, the all-ones model has a smoothing effect on the state at each layer. We note that this generally occurs with any initialization that sets the spectral weights on the same order of magnitude as the affine weights; for instance, the same super-convergence effect occurs when all weights are initialized $U(0,1)$. \mt{We hypothesize that this is because when the spectral weights are of equal magnitude to the affine weights, the function is progressively smoothed as it passes through the model.}

\begin{figure}[t]
    \centering
    \includegraphics[width = 0.35\linewidth]{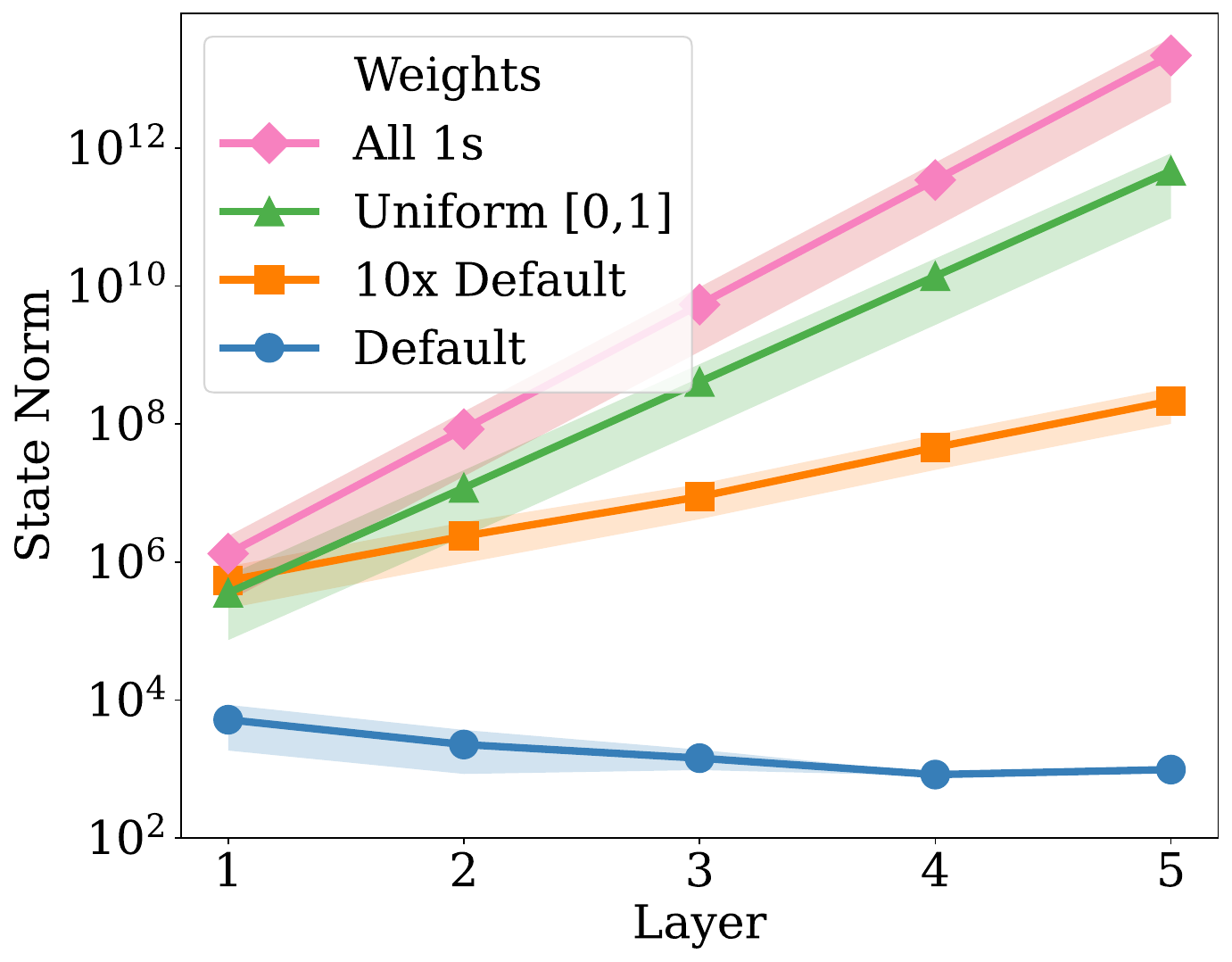}
    \caption{State norm versus layer for various untrained model initializations.}
     \label{fig:state-norm-vs-layer}
\end{figure}

We can also observe the state norm as the state passes through the layers for various settings of the weights. Indeed, for all choices of initialization that we explored except the default setting, the state norm increases exponentially through the layers, while for the default initialization the magnitude stays roughly constant. This phenomenon is illustrated in Figure \ref{fig:state-norm-vs-layer}.

\section{Additional implementation details for error analysis experiments}

All the trained models were trained on an Nvidia P100 GPU for approximately $6$ hours. The evaluation scripts were run on a Mac laptop with an M2 processor. 

\section{Implementation details for adaptive subsampling}
\label{appx:code_details}
Our model has 4 hidden layers, channel width 64 and Fourier cut-off 12. Our results are based on 9000 training samples and 500 test samples. For training with a subsampling scheduler, we include an additional 500 samples for validation. Models are trained for 300 epochs on an Nvidia P100 GPU.

\end{document}